\documentclass[11pt]{amsart}
\usepackage{color,enumitem,graphicx}
\usepackage{amsthm}
\usepackage{amsmath,amssymb, esint}
\usepackage{mathrsfs}
\newtheorem{thm}{Theorem}[section]

\newtheorem{lem}[thm]{Lemma}
\newtheorem{prop}[thm]{Proposition}
\theoremstyle{definition}
\newtheorem{defn}[thm]{Definition}
\theoremstyle{remark}

\numberwithin{equation}{section}

\newcommand{\R}{\mathbb{R}}

\newcommand{\RR}{\mathbb{R}}
\newcommand{\GG}{\mathbb{G}}
\newcommand{\grad}{\nabla_\gamma}
\newcommand{\lap}{\Delta_{\gamma, p}}

\newcommand{\abs}[1]{\left\vert#1\right\vert}

\newcommand{\norm}[1]{\left\Vert#1\right\Vert}

\begin{document}

      \title[]
{Quasilinear problems with critical Sobolev exponent for the
Grushin $p$-Laplace operator}

    \author[S.\ Gandal, A.\ Loiudice, J.\,Tyagi]
     {Somnath Gandal, Annunziata Loiudice$^*$, Jagmohan Tyagi }
    \address{Somnath Gandal \hfill\break
        Central University of Rajasthan, Ajmer \newline NH-8, Bandar Sindri, Dist-Ajmer-305817, Rajasthan, India
        .}
    \email{somnathg@curaj.ac.in, gandalsomnath@gmail.com}
    \address{Annunziata Loiudice \hfill\break
    Universit\`{a} degli Studi di Bari Aldo Moro\newline
    Dipartimento di Matematica, Via Orabona 4, 70125 Bari (Italy)}
    \email{annunziata.loiudice@uniba.it}
    \address{Jagmohan\,Tyagi \hfill\break
       Indian Institute of Technology Gandhinagar \newline
       Palaj, Gandhinagar Gujarat, India-382355.}
   \email{jtyagi@iitgn.ac.in, jtyagi1@gmail.com}

    \begin{abstract}
 We study the
following class of quasilinear degenerate elliptic equations with
critical nonlinearity
\begin{align*}
    \begin{cases}-\Delta_{\gamma,p} u= \lambda |u|^{q-2}u+\abs{u}^{p_{\gamma}^{*}-2}u  & \text{ in } \Omega\subset \RR^N, \\ u=0 & \text{ on } \partial \Omega, \end{cases}
\end{align*}
where $\Delta_{\gamma, p}v:=\sum_{i=1}^N X_i(|\nabla_\gamma
u|^{p-2}X_i u)$ is the Grushin $p$-Laplace operator, $z:=(x, y)
\in \R^N$, $N=m+n,$ $m,n \geq 1,$, where  $\grad=(X_1, \ldots,
X_N)$ is the Grushin gra\-dient, defined as the system of vector
fields $X_i=\frac{\partial}{\partial x_i}, i=1, \ldots, m$,
$X_{m+j}=|x|^\gamma \frac{\partial}{\partial y_j}, j=1, \ldots,
n$, where $\gamma>0$. Here, $\Omega \subset \mathbb{R}^{N}$ is a
smooth bounded domain such that $\Omega\cap \{x=0\}\neq
\emptyset$, $\lambda>0$, $q \in [p,p_\gamma^*)$, where
$p_{\gamma}^{*}=\frac{pN_\gamma}{N_\gamma-p}$ and
$N_\gamma=m+(1+\gamma)n$ denotes the homogeneous dimension
attached to the Grushin gradient. The results extends to the
$p$-case the Brezis-Nirenberg type results in
Alves-Gandal-Loiudice-Tyagi [J. Geom. Anal. 2024, 34(2),52]. The
main crucial step is to preliminarily establish the existence of
the extremals for the involved Sobolev-type inequality
$$
\int_{\RR^N} |\grad u|^p dz \geq S_{\gamma,p} \left ( \int_{\RR^N}
|u|^{p_\gamma^*} dz \right )^{p/p_\gamma^*}
$$
and their qualitative behavior as positive entire solutions to the
limit problem
\begin{equation*}
-\Delta_{\gamma,p} u= u^{p_{\gamma}^{*}-1}\quad \mbox{on}\, \RR^N,
\end{equation*}
whose study has independent interest.

    \end{abstract}
\thanks{A. Loiudice is the corresponding author and she is supported by the INdAM -
GNAMPA Project 2025 "Aspetti qualitativi per equazioni nonlineari
subellittiche su gruppi di Lie" - CUP E5324001950001; she is also
supported by the National Centre on HPC, Big Data and Quantum
Computing, MUR: CN00000013 - CUP: H93C22000450007.}
 \thanks{J. Tyagi thanks DST/SERB for the financial support under the grant CRG/2020/000041.}
    \maketitle

 \tableofcontents
 \section{Introduction}
Our aim in this paper is to study quasilinear Brezis-Nirenberg
type problems related to the so-called $p$-Laplace Grushin
operator. Precisely, we deal with the following nonlinear
subelliptic equations
\begin{align} \label{e:main}
    \begin{cases}-\Delta_{\gamma,p} u=\lambda
    |u|^{q-2}u+\abs{u}^{p_{\gamma}^{*}-2}u
     & \text{ in } \Omega\subset \RR^N =\RR^m\times \RR^n, \\ u=0 & \text{ on } \partial \Omega, \end{cases}
\end{align}
where the Grushin $p$-Laplace operator, $1<p<N_\gamma$, is defined
as
\begin{equation}\label{e:lap}
\Delta_{\gamma, p}u:=\sum_{i=1}^N X_i(|\nabla_\gamma u|^{p-2}X_i
u),
\end{equation}
 being  $\grad=(X_1, \ldots, X_N)$ the usual Grushin gradient,
defined as the system of vector fields
\begin{equation}\label{e:grad}
X_i=\frac{\partial}{\partial x_i}, i=1, \ldots, m, \quad
X_{m+j}=|x|^\gamma \frac{\partial}{\partial y_j}, j=1, \ldots,
n,\quad \gamma>0,
\end{equation}
 and $N_\gamma$ denotes the homogeneous space
dimension naturally attached to the Grushin geometry, i.e.
$$N_\gamma=m+(1+\gamma)n.$$
 Here, $\Omega
\subset \mathbb{R}^{N}$ is a smooth bounded domain intersecting
the degeneration set of the operator, i.e. such that $\Omega \cap
\{x=0\}\neq \emptyset$, the exponent
$$p_{\gamma}^{*}=\frac{pN_\gamma}{N_\gamma-p}$$ is the critical
Sobolev exponent in this context, and $p\leq q<p_\gamma^*$.

This work is the natural continuation of the analysis started by
Alves-Gandal-Loiudice-Tyagi in \cite{AGLT}, where the semilinear
Brezis-Nirenberg type problem for the Grushin operator, that is
the case $p=2$ in \eqref{e:main}, was settled and deeply
investigated, extending to this degenerate-elliptic context the
celebrated results in \cite{Bre}.

In this paper we will explore the quasilinear case $1<p<N_\gamma$,
extending to the Grushin context the seminal results obtained for
the Euclidean $p$-Laplacian operator by Garc\`{i}a Azorero-Peral
in \cite{GP1}, \cite{GP2}.

We emphasize that in the recent literature a great deal of
interest has been devoted to the Grushin operator, especially in
the semilinear case $p=2$ (see e.g. \cite{AFL}, \cite{Alv, Alv
new, AGLT}, \cite{Chen}, \cite{Duo0}, \cite{Duo1}, \cite{G},
\cite{Kog}, \cite{KLV}, \cite{Lamb}, \cite{Liu},
\cite{LoiuNonlinear}, \cite{LoiuAMPA}, \cite{LT}, \cite{MMS1},
\cite{Tri}, \cite{Wang}). Less is known about the general
quasilinear case $1<p<N_\gamma$, especially when the defining
vector fields \eqref{e:grad} are not of H\"{o}rmander's type, that
is, when $\gamma$ is not an even integer. Some results in the
quasilinear case can be found in \cite{B}, \cite{BG}, \cite{Hua
Chen}, \cite{DFZ}, \cite{HY}, \cite{LiuXiao}, \cite{MMS},
\cite{Wei}.

Let us introduce our results. Let us define the space
$\mathring{W}^{1, p}_\gamma(\Omega)$ as the completion of
$C_0^\infty(\Omega)$ with respect to norm $$\|u\|:=\left(
\int_\Omega |\grad u|^p dz \right)^{\frac{1}{p}}.$$ We will look
for weak solutions of equation \eqref{e:main}, i.e. functions
$u\in \mathring{W}^{1, p}_\gamma(\Omega)$ such that
$$
\int_\Omega  |\grad u|^{p-2} \langle \grad u, \grad \phi\rangle \,
dz= \lambda \int_\Omega |u|^{q-2}u \,\phi \,dz + \int_\Omega
|u|^{p^*_\gamma-2} u \,\phi \,dz, \quad \forall \phi \in
\mathring{W}^{1, p}_\gamma(\Omega).
$$
The starting point of the variational formulation of critical
problems of the type \eqref{e:main} is the following $p$-Sobolev
inequality, which follows from the embedding results in \cite{Kog}
(see also \cite[Prop. 17]{LiuXiao}): there exists a constant
$C=C(\gamma,p)>0$ such that
\begin{equation}\label{e:Sobol}
\left ( \int_{\RR^N} |u|^{p^*_\gamma} dz \right
)^{p/p^*_\gamma}\leq C \int_{\RR^N} |\grad u|^{p} \, dz ,\quad
\forall u\in C_0^\infty(\RR^N).
\end{equation}
The first crucial issue that we are lead to investigate in order
to apply the Brezis-Nirenberg approach introduced in \cite{Bre} is
to establish the existence and the qualitative behavior of the
extremal functions in the above inequality. We emphasize that the
explicit form of such extremals is unknown, and only partial
results are available in the case $p=2$ (see \cite{Beck},
\cite{Mon}, \cite{Dou}). We will circumvent this difficulty
arguing as in \cite{Loiu00}.

To this aim, let us introduce the natural \emph{gauge} associated
to the Grushin operator. For $z=(x,y)\in \RR^m\times \RR^n$, let
\begin{equation}\label{defd}
d(z)=(|x|^{2(\gamma+1)}+(\gamma+1)^2|y|^2)^{\frac{1}{2(\gamma+1)}}.
\end{equation}
The role of the function $d$ with respect to the Grushin geometry
is recalled in Section \ref{sec:prelim}. Let us denote by
$\mathcal{D}_\gamma^{1,p}(\RR^N)$ the space defined as the
completion of $C_0^\infty(\RR^N)$ with respect to  $
 \norm{v}_{\mathcal{D}_{\gamma}^{1,p}\left(\mathbb{R}^N\right)}= \left(\int_{\mathbb{R}^{N}}\abs{\nabla_{\gamma} v}^{p} dz \right)^{\frac{1}{p}}.
 $
Our first main result can be stated as follows.

\begin{thm} \label{theo:extrem}The best constant in the $p$-Sobolev inequality \eqref{e:Sobol}, defined as
\begin{align} \label{Sgamma{mn}}
    \mathcal{S}_{\gamma,p}= \inf_{u \in \mathcal{D}_\gamma^{1,p}\left(\mathbb{R}^N
    \right)}
    \frac{\norm{\nabla_{\gamma}
    u}_{L^p(\mathbb{R}^N)}^{p}}{\norm{u}^p_{L^{p_{\gamma}^{*}}(\mathbb{R}^N)}}=\inf_{ u \in \mathcal{D}^{1,p}_{\gamma}\left(\mathbb{R}^N
    \right), \norm{u}_{L^{p_{\gamma}^{*}}(\mathbb{R}^N)}=1}
    \norm{\nabla_{\gamma} u}_{L^p(\mathbb{R}^N)}^{p}
\end{align}
is achieved.

Moreover, any extremal function $u$ for $S_{\gamma,p}$ satisfies
the following estimates
\begin{equation}
\frac{C^{-1}}{1+d(z)^{\frac{N_\gamma-p}{p-1}}}\leq u(z)\leq
\frac{C}{1+d(z)^{\frac{N_\gamma-p}{p-1}}},
\end{equation}
for a suitable constant $C>0$.
\end{thm}
The proof of Theorem \ref{theo:extrem} will follow from the
analysis performed in Section \ref{sec:exist} and
\ref{sec:limitprob} (see, in particular, Theorem
\ref{theo:minimiz} and Theorem \ref{theo:asympt}). Precisely, in
Section \ref{sec:exist} we develop the concentration-compactness
methods that allow us to prove the existence of the $p$-Sobolev
Grushin extremals, extending the results for the case $p=2$ in
\cite{AGLT}. Then, in Section \ref{sec:limitprob}, we obtain the
sharp asymptotic behavior of such extremals as positive entire
solutions of the limit problem
\begin{equation*}
-\Delta_{\gamma,p} u= u^{p_{\gamma}^{*}-1}\quad \mbox{in}\, \RR^N.
\end{equation*}
by means of a refined regularity analysis, following the scheme in
\cite{LoiuOpt}. Concerning the existence of minimizers for the
best constant in \eqref{Sgamma{mn}}, we point out that, when
$\gamma$ is an even positive integer, the result is a special case
of the existence result in \cite{Hua Chen}, where the case of
homogeneous H\"{o}rmander vector fields is addressed.

Now, by means of the fundamental tools in Theorem
\ref{theo:extrem}, we can tackle the problem of the existence of
positive solutions for problem \eqref{e:main} on bounded domains
$\Omega \subset \RR^N$ intersecting the degeneration set of the
operator. We emphasize that the condition $\Omega \cap \{x=0\}\neq
\emptyset$ ensures that the considered problem with exponent
$p_\gamma^*=\frac{pN_\gamma}{N_\gamma-p}$ turns out to be
``critical" in the usual sense of the Sobolev embeddings, since
for these domains the Sobolev embedding
$\mathring{W}_\gamma^{1,p}(\Omega)\hookrightarrow L^p(\Omega)$ is
compact for $p<p_\gamma^*$ (see \cite{Kog}), but fails to be
compact for $p=p_\gamma^*$, due to the action of the rescalings
 defined in \eqref{e:transf}.

Our first result in this direction concerns the problem with
$p$-linear perturbation, that is
\begin{align}\label{e:mainp}
    \begin{cases}-\Delta_{\gamma,p} u= \lambda |u|^{p-2}u+\abs{u}^{p_{\gamma}^{*}-2}u
    & \text{ in } \Omega\subset \RR^N, \\ v=0 & \text{ on } \partial \Omega. \end{cases}
\end{align}
Denote by $\lambda_1=\lambda_1(\Omega)$ the first eigenvalue of
$-\Delta_{\gamma,p}$ on $\mathring{W}^{1,p}_{\gamma}(\Omega)$,
defined by
 $$ \lambda_1:=\inf _{u \in
\mathring{W}^{1,p}_{\gamma}(\Omega)}\frac {\|\nabla_{\gamma,p}
u\|^p_p}{\|u\|^p_p}. $$ Observe that $\lambda_1>0$, due to the
validity of the continuous embedding
$\mathring{W}_\gamma^{1,p}(\Omega)\hookrightarrow L^p(\Omega)$
(see e.g. \cite{Kog}). We prove the following result.

\begin{thm}\label{boundedmainresult1}
Let $\Omega\subset \mathbb{R}^N$ be a smooth bounded domain,
$\Omega \cap \{x=0\}\neq \emptyset$.  Then, \\
i) if $N_\gamma \geq p^2$, problem \eqref{e:mainp} has a
nontrivial solution $u\in \mathring{W}^{1,
p}_\gamma(\Omega)$ for any $0<\lambda <\lambda_{1}$;\\
ii) if $N_\gamma<p^2$, there exists $\lambda_*>0$ such that
\eqref{e:mainp} has a nontrivial solution $u\in \mathring{W}^{1,
p}_\gamma(\Omega)$ for $\lambda \in (\lambda_*, \lambda_1)$.
\end{thm}
To complete the picture of known results about problem
\eqref{e:mainp}, we quote the bifurcation result in \cite{MMS},
where the existence of solutions in a left neighborhood of any
Dirichlet eigenvalue is proved.

 To complete our analysis, we consider the critical
problem with a general subcritical perturbation
\begin{align}\label{e:mainp1}
    \begin{cases}-\Delta_{\gamma,p} u= \lambda |u|^{q-2}u+\abs{u}^{p_{\gamma}^{*}-2}u
        & \text{ in } \Omega\subset \RR^N, \\ v=0 & \text{ on } \partial \Omega, \end{cases}
\end{align}
where $p<q<p_\gamma^*.$ We prove the following result.

\begin{thm}\label{boundedmainresult2}
    Let $\Omega\subset \mathbb{R}^N$ be a smooth bounded domain,
    $\Omega \cap \{x=0\}\neq \emptyset$ and let $p<q<p_\gamma^*$. Then, \\
    i) if $N_\gamma \geq p^2$, problem \eqref{e:mainp1} has a
    nontrivial solution $u\in \mathring{W}^{1,
        p}_\gamma(\Omega)$ for all $\lambda>0$;\\
    ii) if $N_\gamma<p^2$ and $ p_\gamma^*-\frac{p}{p-1}<q<p_\gamma^*$,
    problem \eqref{e:mainp1} has a nontrivial solution for any
    $\lambda>0$; if $p<q\leq p_\gamma^*-\frac{p}{p-1}$, there exists
    $\lambda_*>0$ such that problem \eqref{e:mainp1} has a nontrivial
    solution for $\lambda>\lambda_*$.
\end{thm}
We emphasize that, in the above existence results, the homogeneous
space dimension $N_\gamma$ plays exactly the same role as the
topological dimension $N$ in the Euclidean elliptic setting. The
same phenomenon was observed in the subelliptic context of
sub-Laplacians on Carnot groups (see \cite{LoiuOpt},
\cite{LoiuTrends} for related results in the stratified context).

The paper is organized as follows: in Section \ref{sec:prelim} we
introduce the functional setting of the Grushin $p$-Laplace
operator; in Section \ref{sec:exist} we prove the existence of
$p$-Sobolev minimizers; then, in Section \ref{sec:limitprob} we
obtain the main qualitative properties of such minimizers, in
particular their exact rate of decay at infinity; finally, in
Section \ref{sec:BNexist} we prove the existence results on
bounded domains stated in Theorem \ref{boundedmainresult1} and
\ref{boundedmainresult2}.

\section{The functional setting}\label{sec:prelim}
In this section we recall some basic properties of the considered
operator and the related functional setting. If $ X =
(X_1,\ldots,X_N)$ is the system of Grushin vector fields defined
in \eqref{e:grad}, that is $X_i=\frac{\partial}{\partial x_i},
i=1, \ldots, m$, $X_{m+j}=|x|^\gamma \frac{\partial}{\partial
y_j}, j=1, \ldots, n$, where $\gamma>0$, $n,m\geq 1$, $n+m=N$, the
nonlinear differential operator
\begin{equation*}
\lap u = \sum_{i=1}^m X_i ( |Xu|^{p-2}X_i u)
\end{equation*}
will be referred to as the Grushin $p$-\emph{Laplace} operator.
Note that for any $c>0$, one has $\lap (cu)=c^{p-1}\lap u$ and
furthermore, since the $X_j$'s are homogeneous of degree one with
respect to the dilations
\begin{equation}\label{e:dil}
\delta_{\rho}(x,y)=(\rho x, \rho^{\gamma+1}y),\ \ \forall (x,y)\in
\RR^m\times \RR^n,\ \rho>0
\end{equation}
the operator $\lap $ is homogeneous of degree $p$ with respect to
$\delta_\rho$, namely
\begin{equation*}
\lap (u\circ \delta_\rho)= \rho^{p} \lap u\circ \delta_\rho, \ \
\rho>0.
\end{equation*}
Due to the homogeneity of the operator with respect to the
anisotropic dilations \eqref{e:dil}, the associated
\emph{homogeneous dimension}
$$N_\gamma=m+(\gamma+1)n$$ plays a crucial role in the analysis of
the operator.

The linear operator $\Delta_\gamma=\Delta_{\gamma,2}$,
corresponding to the case $p=2$, is the so-called Grushin (or
Baouendi-Grushin) operator $$\Delta_\gamma=
\Delta_x+|x|^{2\gamma}\Delta_y,$$ which have been deeply studied
in the literature as a model case of a wide class of homogeneous
degenerate-elliptic operators (see \cite{Baou}, \cite{Gr1, Gr2},
\cite{FL1, FL2, FL3}, \cite{G}).

We emphasize that the Grushin vector fields $X_j's$ \eqref{e:grad}
fall into the general class of H\"{o}rmander vector fields only
when $\gamma$ is a positive even integer. In this case a quite
extensive literature about regularity theory and related
Sobolev-type embeddings is available (see e.g.
Capogna-Danielli-Garofalo \cite{CDG1}).

On the other hand, when $\gamma>0$ is any positive real number,
the Grushin operator can be seen as a special case of the general
class of operators introduced and studied by Franchi and
Lanconelli in \cite{FL1, FL2, FL3}, and Kogoj-Lanconelli in
\cite{Kog}.

Concerning the case of quasilinear equations involving possibly
non-smooth Grushin vector fields, some regularity tools such as
Moser-type estimates and Harnack inequality for can be found in
\cite{DFZ}, \cite{GJ}.

For $z=(x,y)\in \RR^m\times \RR^n$, let
$$
d(z):=(|x|^{2(\gamma+1)}+(\gamma+1)^2|y|^2)^{\frac{1}{2(\gamma+1)}}.
$$
The function $d(z)$ is homogeneous of degree one with respect to
the anisotropic dilations \eqref{e:dil}. We will denote by $B_R$
the set
$$
B_R:=\{z\in \RR^N : \, d(z)<R \}.
$$
It holds that $$|B_R|= |B_1| R^{N_\gamma},$$ where $|A|$ denotes
the Lebesgue measure of any measurable set $A\subset \RR^N$.

 The function
$d(z)$ is related to the fundamental solution at the origin of the
Grushin operator. Precisely, if $p>1$ and $\Gamma_p$ is defined,
for $z\neq 0$, as
$$
\Gamma_p(z):= \begin{cases} d(z)^{\frac{p-N_\gamma}{p-1}}\ \ \mbox{if}\ p\neq N_\gamma\\
-\ln [d(z)]\ \ \mbox{if}\ p=N_\gamma
\end{cases}
$$
by direct computation, it is possible to show that $-\lap
\Gamma_p=0$ on $\RR^N\setminus \{0\}$ and that there exists a
constant $C_p\neq0$ such that
$$-\lap \Gamma_p=C_p \delta_0\ \ \mbox{on}\ \RR^N$$
in the weak sense and $C_p>0$ if and only if $N_\gamma\geq p>1$
(see \cite{G} for the case $p=2$; see \cite{B, BG} for the general
$p$-case). The function $\Gamma_p$ will be useful as a comparison
function in the estimate from below for minimizers in Section
\ref{sec:limitprob}.

The starting point of the variational formulation of our problem
is the validity of the global Sobolev-type inequality
\eqref{e:Sobol}. Such inequality can be obtained as a consequence
of the embedding results by Franchi and Lanconelli in \cite{FL3}
(see \cite{LoiuSobol} for the case $p=2$ and \cite{Kog} for the
general $p$-case). We quote \cite{Kog} for a generalization to the
class of $\Delta_\lambda$-operators; finally, we point out the
recent general results in \cite{Hua Chen}, where a large variety
of $p$-Sobolev-type embeddings is investigated in the general
setting of homogeneous H\"{o}rmander vector fields.

\section{Existence of $p$-Sobolev extremals}\label{sec:exist}

Let $1<p<N_\gamma$ and consider the following Sobolev-type
inequality
\begin{equation}\label{e:Sobol1}
 \int_{\RR^N} |\grad u|^p dz  \geq
S_{\gamma,p} \left ( \int_{\RR^N} |u|^{p_\gamma^*} dz \right
)^{p/p_\gamma^*}\quad \forall u\in C_0^\infty(\RR^N).
\end{equation}
Define the best constant
\begin{align} \label{Sgamma}
    \mathcal{S}_{\gamma,p}:= \inf_{ u \in \mathcal{D}_\gamma^{1,p}\left(\mathbb{R}^N \right), \norm{u}_{L^{p_{\gamma}^{*}}(\mathbb{R}^N)}=1}
    \norm{\nabla_{\gamma} u}_{L^p(\mathbb{R}^N)}^{p}>0.
\end{align}
In this section, our goal is to prove that the infimum defined
above is achieved, for any $m,n\geq 1$. To accomplish our goal, we
extend the proof for the case $p=2$ performed in the previous
paper \cite{AGLT}, whose Euclidean outline can be found in
\cite[Section 1.9]{Wil}.

Let $u\in  \mathcal{D}^{1,p}_{\gamma}\left(\mathbb{R}^N \right)$
and $\rho
>0.$ Define
\begin{equation}\label{e:transf}
u^{e,\rho}(z):= \rho^{\frac{N_{\gamma}-p}{p}}u(\rho x,
\rho^{1+\gamma}y + e),
\end{equation}
 where $z=(x,y) \in \mathbb{R}^{m+n}$ and $
e \in \mathbb{R}^n.$ It is easy seen that the following
invariances hold:
$$\norm{\nabla_{\gamma} u^{e, \rho}}_{L^p(\mathbb{R}^N)}= \norm{\nabla_{\gamma}u}_{L^p(\mathbb{R}^N)} \text{ and }
 \norm{u^{e, \rho}}_{L^{p_\gamma^*}(\mathbb{R}^N)}= \norm{u}_{L^{p_\gamma^*}(\mathbb{R}^N)}.$$
Analogously, it can be verified that the equation satisfied, up to
a stretching constant, by the minimizers for $S_{\gamma,p}$, that
is $-\Delta_{\gamma,p}u=|u|^{p^*_\gamma-2}u$ in $\RR^N$, is
invariant under the translations in the $y$ variable and the
rescaling defined in \eqref{e:transf}. Note that the same problem
is not invariant under general translations in $\RR^N$.

Let us recall some preliminary results and definitions from
measure theory \cite{Wil}. Let $\Omega$ be an open subset of
$\mathbb{R}^N$ and define
$$
    K(\Omega):=\{u \in C(\Omega): \operatorname{supp} u \text { is a compact subset of } \Omega\}
$$
and
$$
    B(\Omega):=\left\{u \in C(\Omega):\norm{u}_{L^{\infty}(\Omega)}:=\sup _{z \in \Omega}|u(z)|<\infty\right\}.
$$
We denote by $C_{0}(\Omega)$ the closure of $K(\Omega)$ in
$B(\Omega)$ with respect to the uniform norm. Following the
approach in \cite{Wil1}, we adopt the following definitions.
\begin{defn}
    A finite measure on $\Omega$ is a continuous linear functional on $C_0(\Omega)$. The norm of the finite measure $\mu$ is defined by
    $$
    \norm{\mu}:=\sup _{\substack{u \in C_0(\Omega) \\\norm{u}_{L^{\infty}(\Omega)}=1}}\abs{\langle\mu, u\rangle} .
    $$ \end{defn}
We denote by $$\mathcal{M}(\Omega)- \text{The space of finite
measures on $\Omega.$ }  $$

$$\mathcal{M}^{+}(\Omega)- \text{The space of positive finite measures on $\Omega.$ } $$

We say that a sequence $\mu_{k} \rightharpoonup \mu $ weakly in
$\mathcal{M}(\Omega)$, whenever
$$
\left\langle\mu_k, u\right\rangle \rightarrow\langle\mu, u\rangle,
\quad \forall u \in C_0(\Omega).
$$
\begin{thm}
    (i) Every bounded sequence of finite measures on $\Omega$ contains a weakly convergent subsequence.\\
    (ii) If $\mu_k \rightharpoonup \mu$ in $\mathcal{M}(\Omega),$ then $\left(\mu_k\right)$ is bounded and
    $$
    \norm{\mu} \leq \liminf \norm{\mu_{k}} .
    $$
    (iii) If $\mu \in \mathcal{M}^{+}(\Omega),$ then
  $B(\Omega)\subset L^1(\Omega, \mu)$ and
    $$
    \norm{\mu}=\sup _{\substack{u \in B(\Omega) \\\norm{u}_{L^{\infty}(\Omega)}=1}}| \langle\mu, u\rangle| =\langle\mu, 1\rangle.
    $$
\end{thm}
The proof of the above theorem can be found in \cite{Wil1}, page
206.

\begin{lem}\label{ccp}
    (The Concentration-Compactness Principle) Let $\left\{u_{k}\right\} \subset \mathcal{D}^{1,p}_{\gamma}\left(\mathbb{R}^N \right) $ be a sequence such that
    \begin{align*}
        u_k \rightharpoonup  u & \text { in } \mathcal{D}^{1,p}_{\gamma}\left(\mathbb{R}^N\right), \\
        \abs{\nabla_{\gamma}\left(u_k-u\right)}^p \rightharpoonup \mu & \text { in } \mathcal{M}\left(\mathbb{R}^N\right), \\
        \abs{u_k-u}^{p_{\gamma}^{*}} \rightharpoonup  \nu & \text { in } \mathcal{M}\left(\mathbb{R}^N\right), \\
        u_k \rightarrow u & \text { a.e. on } \mathbb{R}^N
    \end{align*}
    and define
    \begin{align}\label{ccp0}
        \mu_{\infty}:=\lim _{R \rightarrow \infty} \limsup_{k \rightarrow \infty} \int_{d(z) > R}\abs{\nabla_{\gamma} u_k}^p dz, \quad
         \nu_{\infty}:=\lim _{R \rightarrow \infty} \limsup_{k \rightarrow \infty} \int_{d(z)>R}\abs{u_k}^{p_{\gamma}^*} dz,\end{align}
    then it follows that
    \begin{align}
        \norm{\nu}^{\frac{p}{p_{\gamma}^*}} \leq \mathcal{S}_{\gamma,p}^{-1}\,\norm{\mu}, \label{ccp1}\\
        \nu_{\infty}^{\frac{p}{p_{\gamma}^*}} \leq \mathcal{S}_{\gamma,p}^{-1}\,\mu_{\infty}, \label{ccp2}\\
        \limsup_{k \rightarrow \infty}\norm{\nabla_{\gamma} u_k}_{p}^p= \norm{\nabla_{\gamma} u}_{p}^p+\norm{\mu}+\mu_{\infty}, \label{ccp3}\\
        \limsup_{k \rightarrow \infty}\norm{u_k}_{p_{\gamma}^{*}}^{p_{\gamma}^*}=\norm{u}_{p_{\gamma}^*}^{p_{\gamma}^*}+\norm{\nu}+\nu_{\infty}. \label{ccp4}
    \end{align}
    Moreover, if $u=0$ and $\norm{\nu}^{\frac{p}{p_{\gamma}^*}} = \mathcal{S}_{\gamma,p}^{-1}\norm{\mu},$ then the measures $\mu$ and $\nu$ are concentrated at a single point.
\end{lem}
\begin{proof}
    Case (i) Assume first $u=0.$ Let $\phi \in C_{0}^{\infty}(\mathbb{R}^N),$ then from the Sobolev embedding, we have
    $$\mathcal{S}_{\gamma,p} \left(\int_{\mathbb{R}^N} \abs{\phi u_{k}}^{p_{\gamma}^{*}} dz \right)^{\frac{p}{p_{\gamma}^{*}}}
     \leq \int_{\mathbb{R}^{N}} \abs{\nabla_{\gamma} (\phi u_{k})}^p dz .$$

    Now, using the H\"{o}lder inequality and the convergence $u_{k} \rightarrow 0$ in $L^{p}_{loc}(\mathbb{R}^N),$ we obtain
    \begin{align}\label{ccp6}
        \mathcal{S}_{\gamma,p} \left(\int_{\mathbb{R}^{N}} \abs{\phi}^{p_{\gamma}^{*}} d\nu \right)^{\frac{p}{p_{\gamma}^{*}}}
        \leq  \int_{\mathbb{R}^{N}} \abs{\phi}^p d\mu.\end{align}
    Thus \ref{ccp1} follows.\\
    For $R>0,$ let $\phi_{R} \in C^{1}(\mathbb{R}^N)$ be such that $0 \leq \phi_{R} \leq 1 $ on $\mathbb{R}^N,$ $\phi_{R} = 1$ in $B^c_{R+1}:= \mathbb{R}^N\setminus B_{R+1}$ and $\phi_{R}=0$ in $B_{R}.$ Again, using Sobolev inequality, we obtain
    $$\mathcal{S}_{\gamma,p} \left(\int_{\mathbb{R}^N} \abs{\phi_{R} u_{k}}^{p_{\gamma}^{*}} dz \right)^{\frac{p}{p_{\gamma}^{*}}}
    \leq \int_{\mathbb{R}^{N}} \abs{\nabla_{\gamma} (\phi_{R} u_{k})}^p dz .$$
    Now, using  H\"{o}lder's inequality and the convergence $u_{k} \rightarrow 0$ in $L^{p}_{loc}(\mathbb{R}^N),$ we obtain
    \begin{align}\label{ccp5}
        \mathcal{S}_{\gamma,p} \limsup_{k \rightarrow \infty} \left(\int_{\mathbb{R}^N} \abs{\phi_{R} u_{k}}^{p_{\gamma}^{*}} dz
         \right)^{\frac{p}{p_{\gamma}^{*}}}  \leq  \limsup_{k \rightarrow \infty} \int_{\mathbb{R}^{N}} \abs{\nabla_{\gamma}u_{k}}^p
         \abs{\phi_R}^p dz .
    \end{align}
    It is easy to check that
    \begin{align} \label{INE1}
        \int_{B^c_{R+1}} \abs{\nabla_{\gamma} u_{k}}^p dz \leq \int_{\mathbb{R}^{N}} \abs{\nabla_{\gamma} u_{k}}^p \phi_{R}^p dz
        \leq \int_{B^c_{R}} \abs{\nabla_{\gamma} u_{k}}^p dz
    \end{align}
    and
    \begin{align} \label{INE2}
        \int_{B^c_{R+1}} \abs{u_{k}}^{p_{\gamma}^{*}} dz \leq \int_{\mathbb{R}^N} \abs{u_{k}}^{p_{\gamma}^{*}}\phi_{R}^{p_{\gamma}^*} dz
        \leq \int_{B^c_{R}} \abs{u_{k}}^{p_{\gamma}^{*}}dz.
    \end{align}
Now, \eqref{INE1} and \eqref{INE2} combine with (\ref{ccp0}) to
give
    \begin{align*}
        \mu_{\infty}=\lim _{R \rightarrow \infty} \limsup_{k \rightarrow \infty} \int_{B_R^c}\abs{\nabla_{\gamma} u_k}^p \phi_{R}^p dz ,
        \quad \nu_{\infty}=\lim _{R \rightarrow \infty} \limsup_{k \rightarrow \infty} \int_{B_R^c}\abs{u_k}^{p_{\gamma}^*}
        \phi_{R}^{p_{\gamma}^*}dz.
    \end{align*}
    Inequality (\ref{ccp2}) then follows from \ref{ccp5}. \\
    Assume, moreover, that $\norm{\nu}^{\frac{p}{p_{\gamma}^*}} = \mathcal{S}_{\gamma,p}^{-1}\norm{\mu}.$ From (\ref{ccp6}),

    \begin{align*}
        \mathcal{S}^{\frac{1}{p}}_{\gamma,p} \left(\int_{\mathbb{R}^{N}} \abs{\phi}^{p_{\gamma}^{*}} d\nu \right)^{\frac{1}{p_{\gamma}^{*}}}
         \leq  \left(\int_{\mathbb{R}^{N}} \abs{\phi}^pd\mu \right)^{\frac{1}{p}}.
    \end{align*}
    Since $\frac{p}{N_{\gamma}}+ \frac{p}{p_{\gamma}^{*}}=1,$ we use H\"{o}lder's inequality to obtain
    \begin{align*}
        \left(\int_{\mathbb{R}^{N}} \abs{\phi}^p d\mu \right)^{\frac{1}{p}} \leq \norm{\mu}^{\frac{1}{N_{\gamma}}} \left(\int_{\mathbb{R}^{N}}
        \abs{\phi}^{p_{\gamma}^{*}} d\mu \right)^{\frac{1}{p_{\gamma}^{*}}} .
    \end{align*}
    Therefore
    \begin{align*}
        \mathcal{S}^{\frac{1}{p}}_{\gamma, p} \left(\int_{\mathbb{R}^{N}} \abs{\phi}^{p_{\gamma}^{*}} d\nu \right)^{\frac{1}{p_{\gamma}^{*}}}
         \leq \norm{\mu}^{\frac{1}{N_{\gamma}}} \left(\int_{\mathbb{R}^{N}} \abs{\phi}^{p_{\gamma}^{*}} d\mu \right)^{\frac{1}{p_{\gamma}^{*}}},
    \end{align*}
    leading to
    $$
    \nu = \mathcal{S}^{-\frac{p^{*}_{\gamma}}{p}}_{\gamma,p} \norm{\mu}^{\frac{p}{N_{\gamma}-p}}\mu.
    $$
    It follows from (\ref{ccp6}), for each $\phi \in C_{0}^{\infty}(\mathbb{R}^N),$
    \begin{align*}
        \left(\int_{\mathbb{R}^{N}} \abs{\phi}^{p_{\gamma}^{*}} d\nu \right)^{\frac{1}{p_{\gamma}^{*}}}  \norm{\nu}^{\frac{1}{N_{\gamma}}}
         \leq \left(\int_{\mathbb{R}^{N}} \abs{\phi}^p d\nu \right)^{\frac{1}{p}}.
    \end{align*}
    Therefore, for each open set $\Omega \subset \mathbb{R}^N,$
    $$
    \nu (\Omega)^{\frac{1}{p_{\gamma}^*}} \nu(\mathbb{R}^N)^{\frac{1}{N_{\gamma}}} \leq \nu(\Omega)^{\frac{1}{p}},
    $$
    showing that the measure $\nu$ is concentrated at a single point.\\
    Case (ii) Consider, now, the general case. Define
    $$w_{k}:= u_{k}-u.$$
     Then $$w_{k} \rightharpoonup 0\,\, \text{ in } \,\mathcal{D}^{1,p}_{\gamma}(\mathbb{R}^N),$$
    $$ \abs{\nabla_{\gamma} u_{k}}^p \rightharpoonup \mu + \abs{\nabla_{\gamma} u}^p\,\, \text{ in } \,\mathcal{M}(\mathbb{R}^N).$$
    Now, using the Br\'{e}zis-Lieb Lemma ( see, e.g., Lemma 1.32 in \cite{Wil}), we have for any $\phi \in K(\mathbb{R}^N),$
    $$
    \int_{\mathbb{R}^N} \phi \abs{u}^{p_{\gamma}^*}dz=\lim _{k \rightarrow \infty}\left(\int_{\mathbb{R}^N} \phi\abs{u_{k}}^{p_{\gamma}^*}dz-
    \int_{\mathbb{R}^N} \phi \abs{w_{k}}^{p_{\gamma}^*} dz \right)
    $$
    that is,
    $$
    \abs{u_{k}}^{p^*} \rightharpoonup \nu+\abs{u}^{p_{\gamma}^*} \text { in } \mathcal{M}\left(\mathbb{R}^N\right).
    $$
    Hence, inequality (\ref{ccp1}) follows from the corresponding inequality for $\left\{w_{k} \right\}.$
    Since
    $$
    \limsup_{k \rightarrow \infty} \int_{B^c_{R}}\abs{\nabla_{\gamma} w_k}^p dz=\limsup_{k \rightarrow \infty}
    \int_{B^c_{R}}\abs{\nabla_{\gamma} u_k}^pdz-\int_{B^c_{R}}\abs{\nabla_{\gamma} u}^p dz,
    $$
    we obtain
    $$
    \lim _{R \rightarrow \infty} \limsup_{k\rightarrow \infty} \int_{B^c_{R}}\abs{\nabla w_k}^p=\mu_{\infty} .
    $$
    By the Brezis-Lieb lemma,
    $$
    \int_{B^c_{R}}\abs{u}^{p_{\gamma}^*}dz=\lim _{k \rightarrow \infty}\left(\int_{B^c_{R}}\abs{u_k}^{p_{\gamma}^*}dz-
    \int_{B^c_{R}}\abs{w_k}^{p_{\gamma}^*}\right)dz
    $$
    and so
    $$
    \lim _{R \rightarrow \infty} \limsup_{k \rightarrow \infty} \int_{B^c_{R}}\abs{w_{k}}^{p_{\gamma}^*}dz=\nu_{\infty}.
    $$
    Inequality (\ref{ccp2}) then follows from the corresponding inequality for $\left\{w_{k}\right\}.$\\

    For every $R>1$, we have

    \begin{align*}
        \limsup_{k \rightarrow \infty} \int_{\mathbb{R}^N} \abs{\nabla_{\gamma} u_k}^p dz &=\limsup_{k \rightarrow \infty}
        \left(\int_{\mathbb{R}^N} \phi_R\abs{\nabla_{\gamma} u_k}^p dz +\int_{\mathbb{R}^N}\left(1-\phi_R\right)
        \abs{\nabla_{\gamma} u_k}^p dz \right) \\
        &=\limsup_{k \rightarrow \infty} \int_{\mathbb{R}^N} \phi_R \abs{\nabla_{\gamma} u_k}^pdz
        +\int_{\mathbb{R}^N} \left(1-\phi_R\right) d \mu+\int_{\mathbb{R}^N}\left(1-\phi_R\right)\abs{\nabla_{\gamma} u}^pdz .
    \end{align*}

When $R \rightarrow \infty$, we obtain, by Lebesgue theorem,
    $$
    \limsup_{k \rightarrow \infty} \int_{\mathbb{R}^N}\abs{\nabla_{\gamma} u_k}^p dz =\mu_{\infty}+\int_{\mathbb{R}^N} d \mu
    +\int_{\mathbb{R}^N}\abs{\nabla_{\gamma} u}^pdz =\mu_{\infty}+\norm{\mu}+\int_{\mathbb{R}^{N}} \abs{\nabla_{\gamma} u}^pdz.
    $$
    The proof of (\ref{ccp4}) is similar.
\end{proof}
\begin{thm}\label{theo:minimiz}
    Let $\left\{u_{k} \right\} \subset \mathcal{D}_{\gamma}^{1,p}(\mathbb{R}^N)$ be a minimizing sequence such that   \begin{align}
        \norm{u_{k}}_{p^{*}_\gamma}=1, \,\, \norm{\nabla_{\gamma} u_{k}}_{p}^{p} \rightarrow \mathcal{S}_{\gamma,p}.
    \end{align}
    Then, there exists a sequence $\left\{e_{k}, \rho_{k} \right\} \subset \mathbb{R}^{n} \times (0, \infty)$ such
    that $\left\{u^{e_{k}, \rho_{k}}_{k} \right\}$ contains a convergent subsequence. In particular, there exists a minimizer
    for $\mathcal{S}_{\gamma,p}.$
\end{thm}
\begin{proof}
    Up to a subsequence, we can assume that $$u_{k} \rightharpoonup u \text{ in } \mathcal{D}_{\gamma}^{1,p}(\mathbb{R}^N), $$ so that
    $$
    \norm{\nabla_{\gamma}u }_{p}^p \leq \liminf_{k\to \infty} \norm{\nabla_{\gamma}u_{k} }_{p}^p = \mathcal{S}_{\gamma,p}.
     $$
    Thus,  $u$ is a minimizer provided $\norm{u}_{p_{\gamma}^{*}} = 1.$ But we know only that $\norm{u}_{p_{\gamma}^{*}} \leq 1.$
    Therefore, we aim to show that $\norm{u}_{L^{p_{\gamma}^{*}}} = 1.$ Define a concentration function
    $$Q_{k}(\rho):= \sup_{w:=(0,e) \in \mathbb{R}^{N}} \int_{B_{\rho}(w)} \abs{u_{k}}^{p_{\gamma}^{*}} dz,$$
     where $B_\rho (w)= \{z\in \RR^N |\, d(z-w)<\rho\}$,
    $d$ being
    the homogeneous norm defined in \eqref{defd}.

    Note that for each $k,$
    $$\lim_{\rho \rightarrow0^{+} } Q_{k}(\rho)=0, \quad \lim_{\rho \rightarrow \infty} Q_{k}(\rho)=1 .$$ Therefore, there exists $\rho_{k}>0$ such that $Q_{k}(\rho_{k})=\frac{1}{2}.$
    Moreover, since $$ \lim_{\abs{e} \rightarrow \infty } \int_{B_{\rho_{k}}(w)} \abs{u_{k}}^{p_{\gamma}^{*}}dz=0,$$ there exists $e_{k} \in \mathbb{R}^n$ such that
    $$  \int_{B_{\rho_{k}}(w_{k})} \abs{u_{k}}^{p_{\gamma}^{*}}dz=Q_{k}(\rho_{k})=\frac{1}{2}, $$ where
    $w_{k}:=(0,e_{k}) \in \mathbb{R}^{N}.$
    Let us define $$v_k:=u_k^{e_k, \rho_{k}}.$$ Hence
    $$
    \norm{v_{k}}_{p_{\gamma}^{*}}=1, \,\, \norm{\nabla_{\gamma} v_{k}}_{p}^{p} \rightarrow \mathcal{S}_{\gamma, p}
    $$
    and
    \begin{align}\label{lion1}
        \frac{1}{2}=\int_{B_{1}}\abs{v_{k}}^{p_{\gamma}^*}dz=\sup _{w=(0,e) \in \mathbb{R}^N} \int_{B_{1}(w)}\abs{v_{k}}^{p_{\gamma}^*} dz.
    \end{align}

    Since $\left\{v_{k} \right\}$ is bounded sequence in $\mathcal{D}_{\gamma}^{1,p}(\mathbb{R}^N),$ up to a subsequence, we may assume that
    \begin{align*}
        v_k \rightharpoonup  v & \text { in } \mathcal{D}^{1,p}_{\gamma}\left(\mathbb{R}^N\right), \\
        \abs{\nabla_{\gamma}\left(v_k-v\right)}^p \rightharpoonup \mu & \text { in } \mathcal{M}\left(\mathbb{R}^N\right), \\
        \abs{v_k-v}^{p_{\gamma}^{*}} \rightharpoonup  \nu & \text { in } \mathcal{M}\left(\mathbb{R}^N\right),\\
        v_k \rightarrow v & \text { a.e. on } \mathbb{R}^N.
    \end{align*}
Now, by an application of the previous Lemma \ref{ccp},
    \begin{align}\label{lion2}
        \mathcal{S}_{\gamma,p}=
        \lim_{k \rightarrow \infty} \int_{\mathbb{R}^N}\abs{\nabla_{\gamma} v_k}^p dz= \mu_{\infty}+\norm{\mu}+\int_{\mathbb{R}^{N}}
        \abs{\nabla_{\gamma} v}^pdz,
    \end{align}
    \begin{align}\label{lion3}
        1=\lim_{k \rightarrow \infty} \norm{v_k}_{p_{\gamma}^{*}}^{p_{\gamma}^*}=\norm{v}_{p_{\gamma}^*}^{p_{\gamma}^*}+\norm{\nu}+\nu_{\infty},
    \end{align}
    where
    $$\mu_{\infty}:=\lim _{R \rightarrow \infty} \limsup_{k \rightarrow \infty} \int_{B^c_{R}}\abs{\nabla_{\gamma} v_k}^pdz,
    \quad \nu_{\infty}:=\lim _{R \rightarrow \infty} \limsup_{k \rightarrow \infty} \int_{B^c_{R}}\abs{v_k}^{p_{\gamma}^*}dz.$$
    Now, using (\ref{lion2}), (\ref{ccp1}), (\ref{ccp2}) and Sobolev inequality, we get
    \begin{equation}\label{3terms} \left(\left(\norm{v}^{p_{\gamma}^{*}}_{p_{\gamma}^{*}} \right)^{\frac{p}{p_{\gamma}^{*}}}
     + \norm{\nu}^{\frac{p}{p_{\gamma}^*}} +\nu_{\infty}^{\frac{p}{p_{\gamma}^{*}}} \right) \mathcal{S}_{\gamma,p}
     \leq \mathcal{S}_{\gamma,p}.
     \end{equation}
    It follows from (\ref{lion3}) and \eqref{3terms}, taking into account that $p/p_\gamma^*<1$, that $\norm{v}^{p_{\gamma}^{*}}_{p_{\gamma}^{*}},$
     $\norm{\nu},$ $\nu_{\infty}$ are either equal to $0$ or $1.$ By (\ref{lion1}), $\nu_{\infty} \leq \frac{1}{2},$ and
     hence $\nu_{\infty}=0.$ If $\norm{\nu}=1$ then $v=0.$ Then by the previous Lemma \ref{ccp}, the measure $\nu$ is concentrated at a
     single point $z_{0}.$

     Now, we claim that $z_0=\left(0,e\right) \in \mathbb{R}^m \times
\mathbb{R}^{n}.$ Indeed, arguing by contradiction,
      suppose instead that $z_0=(a,e),$ where $a \in \mathbb{R}^m\setminus
      \{0\}$.
      Then there exists $\epsilon_{0}>0$ such that $\mathcal{B}_{2\epsilon_{0}}(z_0) \cap \left\{x=0\right\}=
      \emptyset$, where we are denoting by $\mathcal{B}_{r}(z_0)$
      the Euclidean ball with center at $z_0$ and radius $r$.
       Let $0<\epsilon <\epsilon_{0},$ and consider the ball $\mathcal{B}_{\epsilon}(z_0).$ Since $\nu$ is concentrated at $z_0,$ we have
     \begin{align}\label{z0}
     \lim_{k\rightarrow \infty} \int_{\mathcal{B}_{\epsilon}(z_0)} \abs{v_k}^{p_{\gamma}^{*}} dz =1.
     \end{align}
     Now, we distinguish the cases $1<p<N$ and $p \geq N.$ \\
     \textbf{Case i):} Let $1<p<N.$ Note that the Grushin critical exponent $p^*_\gamma$ satisfies
      $p_{\gamma}^{*}<p^{*}=\frac{Np}{N-p},$ where $p^{*}$ is the critical exponent for the usual Sobolev embedding.
      Applying H\"{o}lder's inequality, we get
     \begin{align*}
    \left( \int_{\mathcal{B}_{\epsilon}(z_0)} \abs{v_k}^{p_{\gamma}^{*}} dz \right)^{\frac{p}{p_{\gamma}^{*}}} \leq
    \left |\mathcal{B}_{\epsilon}(z_0)\right|^{\frac{p(p^{*}-p_\gamma^{*})}{p_{\gamma}^{*}p^{*}}}   \left( \int_{\mathcal{B}_{\epsilon}(z_0)} \abs{v_k}^{p^{*}} dz \right)^{\frac{p}{p^{*}}}.
     \end{align*}
   Now, by applying the Sobolev embedding for the ordinary gradient to the right hand
   side we have
     \begin{align}\label{e:embedd}
     \left( \int_{\mathcal{B}_{\epsilon}(z_0)} \abs{v_k}^{p^{*}} dz \right)^{\frac{p}{p^{*}}}& \leq C \int_{\mathcal{B}_{\epsilon_0}(z_0)}
     \left( \abs{\nabla v_k}^p + \abs{v_k}^{p} \right)
     dz\\
     & \leq C \left ( \int_{\mathcal{B}_{\epsilon_0}(z_0)}
    \abs{\nabla_\gamma v_k}^p  dz+ \left |\mathcal{B}_{\epsilon_0}(z_0)\right|^{\frac{p(p_\gamma^{*}-p)}{p p_\gamma^{*}}}
    \int_{\mathcal{B}_{\epsilon_0}(z_0)} \abs{v_k}^{p^*_\gamma} dz \right )
    \leq
    C,\nonumber
     \end{align}
where the constant $C>0$ does not depend on $\epsilon$ and $k$,
but only on $\epsilon_0$. Here we have used that the norms on the
spaces $\mathring{W}^{1,p}_{\gamma}(\mathcal{B}_{\epsilon_0}(z_0))
$ and $\mathring{W}^{1,p}(\mathcal{B}_{\epsilon_0}(z_0))$ are
equivalent, being $\mathcal{B}_{\epsilon_0}(z_0)$ a bounded set
and far from the degeneration set, and again H\"{o}lder's
inequality. Hence, combining the previous two inequalities, we get
     \begin{align}\label{z01}
       \left( \int_{\mathcal{B}_{\epsilon}(z_0)} \abs{v_k}^{p_{\gamma}^{*}} dz \right)^{\frac{p}{p_{\gamma}^{*}}} \leq C
        \left | \mathcal{B}_{\epsilon}(z_0)\right
        |^{\frac{p(p^*-p_\gamma^{*})}{p_{\gamma}^{*}p^{*}}},
     \end{align}
 where the constant $C>0$ does not depend on $\epsilon$ and $k$.

 \noindent  \textbf{Case ii):} Let $p\geq N.$ We have the following continuous embeddings for the usual gradient.
 When $p=N,$ $$\mathring{W}^{1,p}(\mathcal{B}_{\epsilon}(z_0)) \hookrightarrow L^{q}(\mathcal{B}_{\epsilon}(z_0)) \quad \forall q>1.$$ And, when $p>N,$

   $$\mathring{W}^{1,p}(\mathcal{B}_{\epsilon}(z_0)) \hookrightarrow C^{0,\alpha}(\overline{\mathcal{B}_{\epsilon}(z_0)})
   \hookrightarrow L^{\infty}(\mathcal{B}_{\epsilon}(z_0)).$$

\noindent Choose $q>p_{\gamma}^{*},$ then again by H\"{o}lder's
inequality and the above Sobolev embedding, reasoning as in Case
i), we get
 \begin{align}\label{z02}
    \left( \int_{\mathcal{B}_{\epsilon}(z_0)} \abs{v_k}^{p_{\gamma}^{*}} dz \right)^{\frac{p}{p_{\gamma}^{*}}}
    \leq C \left
    |\mathcal{B}_{\epsilon}(z_0)\right|^{\frac{p(q-p_{\gamma}^{*})}{p_{\gamma}^{*}q}},
\end{align}
with $C$ independent on $\epsilon$ and $k$.

Then choosing $\epsilon$ small enough in \eqref{z01} and
\eqref{z02}, the integral in the left hand side can be made
arbitrarily small, independently of $k$, which
 contradicts \eqref{z0}. Hence $z_0=(0,e)$, for some $e\in
 \RR^n$.

 Furthermore,
from (\ref{lion1}) we deduce $$\frac{1}{2}=\sup_{w=(0, e) \in
\mathbb{R}^N} \int_{B_{1}(w)}
     \abs{v_{k}}^{p_{\gamma}^{*}}dz \geq \int_{B_{1}(z_{0})} \abs{v_{k}}^{p_{\gamma}^*} dz \rightarrow \norm{\nu}=1,$$
    a contradiction. Thus, $\norm{v}_{p_{\gamma}^{*}}=1$ and so
    \begin{align*}
       \|\nabla_\gamma v\|_p^p = \mathcal{S}_{\gamma,p}= \lim_{k\to \infty} \norm{\nabla_{\gamma} v_{k}}_{p}^{p}.
    \end{align*}
\end{proof}

Finally,  if $\Omega$ is any open subset of $\RR^N$, let us denote
by $\mathcal{S}_{\gamma,p}(\Omega)$ the best constant of the
embedding $\mathring{H}^{1,p}_\gamma(\Omega) \hookrightarrow
L^{p_{\gamma}^*}(\Omega)$, that is
\begin{align} \label{SOmega}
    \mathcal{S}_{\gamma,p}(\Omega):= \inf_{ u \in \mathring{H}^{1,p}_\gamma(\Omega), \norm{u}_{L^{p_{\gamma}^{*}}(\Omega)}=1}
    \norm{\nabla_{\gamma} u}_{L^p(\Omega)}^{p}>0.
\end{align}
The following result holds.
\begin{prop} \label{PropS0} Let $\Omega$ be a bounded domain such that $\Omega \cap \{x=0\} \not= \emptyset$.
Then, $\mathcal{S}_{\gamma,p}(\Omega)=\mathcal{S}_{\gamma,p}$.
\end{prop}
\begin{proof} By definition of $\mathcal{S}_{\gamma,p}(\Omega)$ given in \eqref{SOmega}, we know that
$\mathcal{S}_{\gamma,p} \leq \mathcal{S}_{\gamma,p}(\Omega)$.
 Now, let $\{u_k\} \subset \mathcal{D}_{\gamma}^{1,p}(\R^N)$ be a minimizing sequence for $\mathcal{S}_{\gamma,p}$, that is,
$$
\|u_k\|_{p_{\gamma}^{*}}=1, \quad \lim_{k \to
\infty}\norm{\nabla_{\gamma} u_{k}}_{p}=\mathcal{S}_{\gamma,p}.
$$
By density, we can assume $u_k\in C_0^\infty(\RR^N)$. Since
$\Omega \cap \{x=0\} \not= \emptyset$, without loss of generality,
we can assume that $0 \in \Omega$. Thereby, there is a sequence
$\rho_k>0$ such that $u_k^{0,\rho_k} \in C_0^\infty(\Omega)$ for
$k$ sufficiently large. Hence,
$$
\mathcal{S}_{\gamma,p}(\Omega) \leq \liminf_{k\to \infty}
\norm{\nabla_{\gamma} u^{0,\rho_k}_{k}}_{p}^{p}=\liminf_{k\to
\infty}\norm{\nabla_{\gamma}
u_{k}}_{p}^{p}=\mathcal{S}_{\gamma,p},
$$
showing that $\mathcal{S}_{\gamma,p}(\Omega) \leq
\mathcal{S}_{\gamma,p}$, and so,
$\mathcal{S}_{\gamma,p}(\Omega)=\mathcal{S}_{\gamma,p}$.
\end{proof}

\section{Qualitative behavior of extremals}\label{sec:limitprob}

In this section, we determine the qualitative behavior of the
extremal functions for the $p$-Sobolev inequality \eqref{e:Sobol1}
as solutions, up to a multiplicative constant, of the associated
limit problem on $\RR^N$. To this aim we adapt the arguments
developed by the second author in \cite{LoiuOpt} for the
$p$-sub-Laplacian case, which in turn borrow ideas from the
Euclidean results by V\'{e}tois \cite{Vetois} and Xiang
\cite{Xiang}.

We recall that the rate of decay of Grushin Sobolev minimizers in
the case $p=2$ was obtained by Loiudice in \cite{LoiuNonlinear} by
exploiting the conformal properties of the Grushin Kelvin-type
transform introduced by Monti and Morbidelli in \cite{Mon1}. In
particular, it was proved in \cite{LoiuNonlinear} that any weak
positive solution $u\in \mathcal{D}_\gamma^{1,2}(\RR^N)$ to
$-\Delta_{\gamma} u=u^{2^*_\gamma-1}$ is bounded and satisfies the
estimate $u(z)\sim d(z)^{2-N_\gamma}$ as $d(z)\to \infty$.
Subsequently, the exact rate of decay of solutions of an analogous
Grushin equation involving a Hardy-type term has been obtained in
\cite{LoiuAMPA}, again using the Kelvin transform. In the present
quasilinear context, as it happens in the Euclidean quasilinear
case, we need different tools and a refined regularity analysis
will be performed to this aim, following \cite{LoiuOpt}. Our
result can be stated as follows.

\begin{thm}\label{theo:asympt} If $u\in
\mathcal{D}_\gamma^{1,p}(\mathbb{R}^N)$ is a nonnegative solution
of
\begin{equation}\label{e:entire}
-\Delta_{\gamma,p} u= u^{p_{\gamma}^{*}-1}\quad in\, \RR^N,
\end{equation}
then $u\in L^\infty(\RR^N)\cap L^{p^*_\gamma/p', \infty}(\RR^N)$
and it satisfies
$$
u(z)\sim d(z)^\frac{p-N_\gamma}{p-1}\quad as \,\, d(z)\to \infty,
$$
where
$d(z)=(|x|^{2(\gamma+1)}+(\gamma+1)^2|y|^2)^{\frac{1}{2(\gamma+1)}}$
is the homogeneous norm defined in \eqref{defd}.
\end{thm}


The proof will be divided in three steps. We begin by determining
the sharp weak Lebesgue regularity of solutions to problem
\eqref{e:entire}, following the approach in \cite{Vetois}.

\subsection{$L^q$-weak regularity}
For the reader's convenience, we briefly recall the definition of
weak Lebesgue spaces, referring to \cite{Graf} for a detailed
introduction. For any $s\in (0, \infty)$ and any open set $\Omega
\subset \RR^N$, we define the space $L^{s, \infty}(\Omega)$ as the
set of all measurable functions $u:\Omega \to \RR$ such that
\begin{equation*}
[u]_{L^{s, \infty}(\Omega)}:= \sup_{t>0} t\cdot
\mu(\{|u|>t\})^{1/s}< \infty,
\end{equation*}
where $\mu(\{|u|>t\})$ denotes the Lebesgue measure of the set
$\{z\in \Omega: |u(z)|>t\}$. The map $[u]_{L^{s, \infty}(\Omega)}$
is a quasi-norm on $L^{s, \infty}(\Omega)$. The optimal regularity
of solutions to equations of the type \eqref{e:entire} in the
scale of weak Lebesgue spaces is obtained in the following
proposition.

\begin{prop}\label{prop:weak}
Let $f: \RR^N \times \RR \to \RR$ be a Carath\'{e}odory function
such that
\begin{equation}\label{e:condf}
|f(z, s)|\leq \Lambda |s|^{p_\gamma^*-1}, \quad \mbox{for all}\
s\in \RR \ \mbox{and a.e.}\ z \in \RR^N.
\end{equation}
Then, any solution $u\in \mathcal{D}_\gamma^{1,p}(\RR^N)$ of the
equation
\begin{equation}\label{e:f}
-\Delta_{\gamma,p} u= f(z,u)\quad in\, \RR^N
\end{equation}
belongs to $L^{q_0, \infty}(\RR^N)\cap L^{\infty}(\RR^N)$, where
$q_0=\dfrac{p_\gamma^*}{p'}$, $p'$ being the conjugate exponent of
$p$.
\end{prop}

\begin{proof}

We first observe that the global boundedness of $u$ can be proved
by Moser's iteration technique applied to the Grushin context,
starting from the $L^{p^*_\gamma}$-summability of $u$. The proof
follows the standard scheme, so we omit the details referring, for
instance, to the Euclidean outline in \cite{V}.

Let us now prove that $u$ belongs to $L^{q_0, \infty}(\RR^N)$,
with $q_0=\dfrac{p_\gamma^*}{p'}$. Let us consider the function
\begin{equation*}
T_t(u):=sgn(u)\cdot \min(|u|, t), \quad t>0.
\end{equation*}
By Sobolev inequality \eqref{e:Sobol1}, we get that
\begin{equation}\label{e:estSobol}
t^{p_\gamma^*}\mu(\{|u|>h\}) \leq \int_{\RR^N} |T_t
(u)|^{p_\gamma^*}\, {\rm d}z\leq S_{\gamma,p}^{-p_\gamma^*/p}
\left ( \int_{|u|\leq t} |\grad u|^p \, {\rm d}z \right
)^{\frac{N_\gamma}{N_\gamma-p}}.
\end{equation}
 On the other hand, by testing equation
\eqref{e:entire} with $T_t(u)$ and using assumption
\eqref{e:condf} on $f$ we get
\begin{equation}
\begin{split}\label{e:test0}
\int_{|u|\leq t} |\grad u|^p {\rm d}z & = \int_{|u|\leq t}
f(\xi,u)\cdot u\, {\rm d}z + t \int_{|u|>t} f(z, u)\cdot
sgn(u)\, {\rm d}z\\
 &\leq \Lambda \left ( \int_{|u|\leq t} |u|^{p_\gamma^*}\, {\rm d}z +  t \int_{|u|>t}
|u|^{p_\gamma^*-1}\, {\rm d}z \right ).
\end{split}
\end{equation}
Now, concerning the terms in the right hand side of
\eqref{e:test0}, we have that
\begin{equation}\label{e:est1}
\int_{|u|\leq t} |u|^{p_\gamma^*}\, {\rm d}z = \int_{\RR^N}
|T_t(u)|^{p_\gamma^*}\, {\rm d}z -t^{p_\gamma^*} \mu(\{|u|>t\})
\end{equation}
and
\begin{equation}
\begin{split}\label{e:est2}
\int_{|u|>t} |u|^{p_\gamma^*-1}\, {\rm d}z &= (p_\gamma^*-1)
\int_0^\infty
s^{p_\gamma^*-2}\mu(\{|u|>\max(s,t)\})\, {\rm d}s\\
&= h^{p_\gamma^*-1}\mu(\{|u|>t\})+ (p_\gamma^*-1) \int_t^\infty
s^{p_\gamma^*-2}\mu(\{|u|>s\}) \, {\rm d}s.
\end{split}
\end{equation}
Hence, from \eqref{e:test0}, \eqref{e:est1} and \eqref{e:est2}, we
get
\begin{equation}\label{e:preSobol}
\int_{|u|\leq t} |\grad u|^p \,{\rm d}z \leq \Lambda \left (
\int_{\RR^N} |T_t (u)|^{p_\gamma^*}\,{\rm d}\xi + (p_\gamma^*-1) t
\int_t^\infty s^{p_\gamma^*-2} \mu (\{|u|>s\}) \, {\rm d}s \right
).
\end{equation}
Then, by \eqref{e:estSobol} and \eqref{e:preSobol}, and taking
into account that $\int_{\RR^N} |T_t(u)|^{p_\gamma^*}\,{\rm d}z
=o(1)$ as $t\to 0$, we get the key estimate
\begin{equation}
\begin{split}\label{e:estLp*}
t^{p_\gamma^*}\mu(\{|u|>t\}) & \leq C \left ( t \int_t^{\infty}
s^{p_\gamma^*-2} \mu (\{|u|>s\}) \, {\rm d}s \right
)^{\frac{N_\gamma}{N_\gamma-p}},
\end{split}
\end{equation}
for small $t>0$, for some constant $C=C(N_\gamma,p,\Lambda)$. Now,
in order to estimate the right hand side of \eqref{e:estLp*} and
simplify notation, let us define
\begin{equation}\label{e:defF}
F(t):=\left ( \int_t^{\infty} f(s)\,{\rm d}s \right
)^{-{\frac{p}{N_\gamma-p}}},\quad \mbox{where}\ \
f(s):=s^{p_\gamma^*-2} \mu (\{|u|>s\}).
\end{equation}
Taking into account definition \eqref{e:defF} and the fact that
$p_\gamma^*-\frac{N_\gamma}{N_\gamma-p}=
\frac{N_\gamma(p-1)}{N_\gamma-p}=\frac{p_\gamma^*}{p'}$, the
estimate \eqref{e:estLp*} can be rewritten as follows
\begin{equation*}
t^{p_\gamma^*/p'}\mu(\{|u|>t\})\leq C F(t)^{-N_\gamma/p}, \
\mbox{for small}\ t>0.
\end{equation*}
It is not difficult to verify that $F$ is a non-decreasing
function and that $F(0)>0$. So, we can conclude that
\begin{equation*}
t^{p_\gamma^*/p'}\mu(\{|u|>t\}) \leq C F(0)^{-N_\gamma/p}
\end{equation*}
for small $t$, which implies, taking also into account the
boundedness of $u$, that $[u]_{L^{p_\gamma^*/p', \infty}(\RR^N)}<
\infty$.
\end{proof}

\medskip

\subsection{Estimate from above} The upper bound estimate in Theorem \ref{theo:asympt}
will follow as a direct consequence of Theorem \ref{theo:fundam}
below. The first step of the proof is the following preliminary
reverse H\"{o}lder's inequality inspired to Xiang \cite[Lemma
2.3]{Xiang}. In what follows, denoted by $B_R$ the $d$-ball with
center at 0 and radius $R$, we define the annuli:
\begin{equation}\label{e:annuli1}
A_R = B_{5R} \setminus \overline{B}_{2R}\quad \mbox{and}\quad
\widetilde{A}_R = B_{6R} \setminus \overline{B}_{R}, \quad R>0.
\end{equation}
 The following uniform
estimate with respect to $R$ holds.

\begin{lem}\label{lemma1}
Let $V\in L^{N_\gamma/p}(\RR^N)$ and let $u\in
\mathcal{D}_\gamma^{1,p}(\RR^N)$ be a nonnegative solution to
\begin{equation} \label{e:dis}
-\lap u \leq V u^{p-1}\quad {\mbox in}\,\, \RR^N.
\end{equation}
Let $t>p_\gamma^*$. Then, there exists $R_0>0$ depending on $t$
such that for any $R\geq R_0$, it holds
\begin{equation}\label{e:estAr}
\left ( \fint_{A_R} u^t \right )^{1/t}\leq C \left (
\fint_{\widetilde{A}_R} u^{p_\gamma^*} \right )^{1/p_\gamma^*},
\end{equation}
 where 
  $C$ is a positive constant depending on $t$, but not
on $R$.
\end{lem}

\begin{proof} Adapting the proof in \cite[Lemma 2.3]{Xiang}, for any $R>0$
and $z\in \RR^N$, we define the rescaled function
\begin{equation*}
v(z):= u(\delta_R z).
\end{equation*}
Then, by \eqref{e:dis} and being $\lap (u \circ \delta_R )(z)=R^p
(\lap u)(\delta_R z)$, $v$ satisfies
\begin{equation} \label{e:disVR}
 -\lap v \leq V_R
v^{p-1} \quad \mbox{in}\,\,\RR^N,
\end{equation}
where $$V_R(z)=R^p V(\delta_R z),$$ for any $z\in \RR^N$. We shall
prove estimate \eqref{e:estAr} for $v$ on $\widetilde{A}_1$.

 Let $v_m=\min (v, m)$, for
$m\geq 1$. For any $\eta \in C_0^\infty (\widetilde{A}_1)$,
$\eta\geq 0$ and $s\geq 1$, the test function $\varphi=\eta^p
v_m^{p(s-1)}v$ into \eqref{e:disVR} gives
\begin{equation}\label{e:test}
\int_{\widetilde{A}_1} |\grad v|^{p-2} \grad v \cdot \grad
\varphi\,\leq \int_{\widetilde{A}_1} V_R v^{p-1}\, \varphi.
\end{equation}
Concerning the left hand side of \eqref{e:test}, it is easy to see
that for any sufficiently small $\delta>0$, there exists
$C_\delta>0$ such that
\begin{equation}\label{e:delta}
\int_{\widetilde{A}_1} |\grad v|^{p-2} \grad v\cdot \grad
\varphi\, \geq (1-\delta)
\frac{p(s-1)+1}{s^p}\int_{\widetilde{A}_1} |\grad(\eta
v_m^{s-1}v)|^p - C_\delta \int_{\widetilde{A}_1} |\grad \eta|^p
v_m^{p(s-1)}v^p.
\end{equation}
So, by choosing $\delta=1/2$ in \eqref{e:delta} and using Sobolev
inequality \eqref{e:Sobol1}, we obtain
\begin{equation}\label{e:Lhs}
\int_{\widetilde{A}_1} |\grad v|^{p-2} \grad v\cdot \grad
\varphi\, \geq C_1 \left ( \int_{\widetilde{A}_1} |\eta
v_m^{s-1}v|^{p\chi} \right )^{1/\chi}-C_2
\int_{\widetilde{A}_1}|\grad \eta|^p v_m^{p(s-1)}v^p,
\end{equation}
for some constants $C_1, C_2>0$ depending on $N_\gamma, p, s$,
where $\chi=p_\gamma^*/p$. On the other hand, by H\"{o}lder's
inequality
\begin{equation}\label{e:Rhs}
\begin{split}
\int_{\widetilde{A}_1} V_R v^{p-1} \varphi \leq
\|V_R\|_{\frac{N_\gamma}{p}, \widetilde{A}_1} \left (
\int_{\widetilde{A}_1} |\eta v_m^{s-1}v|^{p\chi} \right )^{1/\chi}
=  \|V\|_{\frac{N_\gamma}{p}, \widetilde{A}_R } \left (
\int_{\widetilde{A}_1} |\eta v_m^{s-1}v|^{p\chi} \right)^{1/\chi}.
\end{split}
\end{equation}
So, by \eqref{e:test}, \eqref{e:Lhs} and \eqref{e:Rhs}, we get
\begin{equation}\label{e:24}
\left ( \int_{\widetilde{A}_1} |\eta v_m^{s-1}v|^{p\chi} \right
)^{1/\chi} \leq C_3 \int_{\widetilde{A}_1}|\grad \eta|^p
v_m^{p(s-1)}v^p + C_3 \|V\|_{\frac{N_\gamma}{p}, \widetilde{A}_R }
\left ( \int_{\widetilde{A}_1} |\eta v_m^{s-1}v|^{p\chi} \right
)^{1/\chi}
\end{equation}
for some constant $C_3=C_3(Q, p, s)>0$.

Now, fix $t>p_\gamma^*$ and let $k\in \mathbb{N}$ such that
$p\chi^k \leq t \leq p \chi^{k+1}$. Then, there exists a positive
constant $C_3=C_3(N_\gamma, p, t)$ such that \eqref{e:24} holds
for all $1\leq s\leq \chi^k$.

Since $V\in L^{N_\gamma/p}(\RR^N)$, there exists $R_0>0$ such that
\begin{equation}\label{e:C3}
C_3 \|V\|_{\frac{N_\gamma}{p}, \widetilde{A}_R } \leq 1/2 \quad
\mbox{ for\,any}\,\, R\geq R_0.
\end{equation}
Therefore, for all $R\geq R_0$, it holds
\begin{equation*}
\left ( \int_{\widetilde{A}_1} |\eta v_m^{s-1}v|^{p\chi} \right
)^{1/\chi} \leq C \int_{\widetilde{A}_1}|\grad \eta|^p
v_m^{p(s-1)}v^p,
\end{equation*}
for all $1\leq s\leq \chi^k$, where $C>0$ depends only on
$N_\gamma, p, t$.

Now, by choosing an appropriate cut-off function $\eta$ and
applying Moser's iteration technique, after finitely many
iterations we can conclude that
\begin{equation}\label{e:equav}
\left ( \int_{A_1} v^t \right )^{1/t}\leq C \left (
\int_{\widetilde{A}_1} v^{p_\gamma^*} \right )^{1/p_\gamma^*}
\end{equation}
for $R\geq R_0$, where we recall that $v(z)= u (\delta_R z)$ and
$C$ does not depend on $R$. By a simple change of variable,
\eqref{e:estAr} follows from \eqref{e:equav}.
 \end{proof}

\bigskip

The next crucial step is an estimate of the $L^\infty$-norm of the
solutions to \eqref{e:entire} on annuli by means of the
 sharp $L^{q_0}$-weak norm on larger annuli, which in fact
 provides
 the sharp rate of decay of solutions at infinity.

 In what follows, we shall indicate by
\begin{equation}\label{e:annuli2}
 \widehat{A}_R=B_{4R}\setminus
 \overline{B}_{3R}
 ,\quad  R>0,
 \end{equation}
and $A_R$ will denote, as before, the larger annulus $B_{5R}
\setminus \overline{B}_{2R}$.
\begin{thm}\label{theo:fundam}
Let $u\in \mathcal{D}_\gamma^{1,p}(\RR^N)$ be a solution of
\eqref{e:f} under the assumption \eqref{e:condf}. Let
$q_0=\frac{p_\gamma^*}{p'}$ be the sharp $L^q$-weak summability
exponent found in Prop. \ref{prop:weak}. Then, there exist
constants $R_0, C>0$, such that for any $R\geq R_0$
\begin{equation}\label{e:sharp}
\sup_{\widehat{A}_R} |u| \leq
\dfrac{C}{|A_R|^{\frac{1}{q_0}}}\,[u]_{L^{q_0, \infty}(A_R)}
\end{equation}
where $C$ does not depend on $R$.
\end{thm}

\begin{proof} First, notice that, if $u\in
\mathcal{D}_\gamma^{1,p}(\RR^N)$ is a solution of equation
\eqref{e:f} under the assumption \eqref{e:condf}, then by Kato's
inequality adapted to the subelliptic Grushin context (see e.g.
\cite{DM}), $|u|$ satisfies
\begin{equation*} -\lap |u| \leq |f(z,u)| \leq \Lambda  |u|^{p_\gamma^*-1} \quad
\mbox{in}\,\, \RR^N,
\end{equation*}
from which
\begin{equation}\label{e:V}
-\lap |u| \leq  V |u|^{p-1},\qquad \mbox{where}\,\, V=\Lambda
|u|^{p_\gamma^*-p}.
\end{equation}
Notice that $V\in L^{q}(\RR^N)$ for any $q\geq N_\gamma/p$, being
$u\in L^t(\RR^N)$ for any $t\geq p_\gamma^*$.

Let us set, as before,
\begin{equation*}
v(z):= |u(\delta_R z)|,\quad \,\, R>0,\, z\in \RR^N.
\end{equation*}
Then, in particular, $v$ weakly satisfies the inequality
\begin{equation} \label{e:VR}
 -\lap v \leq V_R
v^{p-1} \quad \mbox{in}\,\, A_1,
\end{equation}
where $V_R(z)=R^p \,V(\delta_R z)$, with $V$ as in \eqref{e:V}.

Let $t> p^*$ be fixed. Thus, in particular, $ V_R \in
L^{t_0}(\RR^N)$, for
$t_0=\frac{t}{p_\gamma^*-p}>\frac{N_\gamma}{p}$. So, by
subelliptic Moser-type estimates for quasilinear
Schr\"{o}dinger-type equations of type \eqref{e:VR} (see e.g.
\cite{CDG1} for H\"{o}rmander-type operators, Di Fazio et al.
\cite{DFZ} for general Grushin-type vector fields), we have that,
for any $q>0$ the following estimate holds for $v_R$
\begin{equation}\label{e:Cap}
\sup_{\widehat{A}_1} v \leq C \left ( \fint_{A_1} v^{q} \right
)^{1/q},
\end{equation}
where $C=C(N_\gamma, q, \|V_R\|_{L^{t_0}(A_1)})$.

The crucial observation, here, is that the norm
$\|V_R\|_{L^{t_0}(A_1)}$ is uniformly bounded with respect to $R$,
for sufficiently large $R$. Precisely, if we choose $R_0>0$ so
that \eqref{e:C3} holds for $V=\Lambda |u|^{p_\gamma^*-p}$, there
exists a constant $C>0$ depending on $N_\gamma,p,t_0, \Lambda$
such that
\begin{equation}\label{e:estVR}
\|V_R\|_{L^{t_0}(A_1)} \leq
C\|u\|^{p_\gamma^*-p}_{L^{p_\gamma^*}(\RR^N)}\quad \forall R\geq
R_0.
\end{equation}
Indeed, by the definition of $V_R$ and by Lemma \ref{lemma1}
applied to \eqref{e:V}, we get that, for any $R\geq R_0$
\begin{equation*}
\begin{split}
\|V_R\|_{L^{t_0}(A_1)}&=R^{p-\frac{N_\gamma}{t_0}}\|V\|_{L^{t_0}(A_R)}\\
& =\Lambda \,R^{p-\frac{N_\gamma}{t_0}}\|u\|_{L^{t}(A_R)}^{p_\gamma^*-p}\\
& \leq C
R^{p-\frac{N_\gamma}{t_0}-(\frac{N_\gamma}{p_\gamma^*}-\frac{N_\gamma}{t})(p_\gamma^*-p)}\|u\|^{p_\gamma^*-p}_{L^{p_\gamma^*}(\widetilde{A}_R)}\\
& \leq C \|u\|^{p_\gamma^*-p}_{L^{p_\gamma^*}(\RR^N)},
\end{split}
\end{equation*}
with $C>0$ not depending on $R$, where we have used that
$p-\frac{N_\gamma}{t_0}-(\frac{N_\gamma}{p_\gamma^*}-\frac{N_\gamma}{t})(p_\gamma^*-p)=0$.
Therefore, the constant $C$ in \eqref{e:Cap} can be made
independent on $R$, for $R\geq R_0$.

Now, equation \eqref{e:Cap} is equivalent, by rescaling, to
\begin{equation}\label{e:cover}
\sup_{\widehat{A}_R} |u| \leq C \left ( \fint_{A_R} |u|^{q} \right
)^{1/q},
\end{equation}
where $C$ depends on $q$, but not on $R$, for sufficiently large
$R$.

Now, choosing $q$ in \eqref{e:cover} so that
$0<q<q_0=p_\gamma^*/p'$, by H\"{o}lder's inequality for weak
Lebesgue norms (see Grafakos \cite{Graf}) we have
\begin{equation}\label{e:holder}
\left (\int_{A_R} |u|^{q} \right )^{1/q}\leq C_{q,q_0}
|A_R|^{1/q-1/q_0}[u]_{L^{q_0, \infty}(A_R)}.
\end{equation}
Henceforth, by \eqref{e:cover} and \eqref{e:holder}, estimate
\eqref{e:sharp} follows.
\end{proof}

\medskip

We are now able to prove the estimate from above in Theorem
\ref{theo:asympt}.
\medskip\\
\noindent \textbf{Proof of Theorem \ref{theo:asympt} - estimate
from above}. From Theorem \ref{theo:fundam}, by taking into
account that $|A_R|\sim R^{N_\gamma}$ and that, by Proposition
\ref{prop:weak}, $u\in L^{q_0, \infty}(\RR^N)$, the estimate from
above  follows by letting $R=\frac{2}{7}d(z)$ in \eqref{e:sharp},
for $d(z)\geq \frac{7}{2} \,R_0$.
\bigskip

\subsection{Estimate from below} To prove the estimate from below in Theorem
\ref{theo:asympt}, we first introduce a comparison principle on
exterior domains for the Grushin $p$-Laplace operator, adapting
the proof for the ordinary $p$-Laplacian in Xiang \cite[Theorem
3.3]{Xiang}. 
We start with a technical lemma.
\begin{lem}\label{lemma:techn}
Let $u,v$ any two positive and weakly differentiable functions
with respect to the Grushin gradient on a domain $\Omega$.
The following inequalities hold:\\
i) if $p\geq 2$, then\\
$$|\grad u|^{p-2}\grad u\cdot \grad \left ( u-\frac{v^p}{u^p}u
\right ) + |\grad v|^{p-2}\grad v\cdot \grad \left
(v-\frac{u^p}{v^p}v \right ) \geq C_p (u^p+v^p) |\grad \log u
-\grad \log v|^p;$$\\
ii) if $1<p<2$, then\\
\begin{equation*}
\begin{split}
 |\grad u|^{p-2}\grad u\cdot &\grad \left ( u-\frac{v^p}{u^p}u
\right ) + |\grad v|^{p-2}\grad v\cdot \grad \left
(v-\frac{u^p}{v^p}v \right ) \\
& \geq C_p (u^p+v^p) \frac{|\grad \log u -\grad \log v|^2}{
(|\grad \log u |+ |\grad \log v|)^{2-p}}
\end{split}
\end{equation*}
where $C_p$ are positive constants depending only on $p$.
\end{lem}
\begin{proof}
The proof straightforwardly follows the proof for the ordinary
gradient in \cite[Lemma 3.1]{Xiang} and it will be omitted.
\end{proof}
 \begin{thm}\label{theo:comparison}
Let $\Omega\subset \GG$ be an exterior domain and let $f \in
L^{N_\gamma/p}(\Omega)$. Let $v\in
\mathcal{D}_\gamma^{1,p}(\Omega)$ be a weak subsolution of $-\lap
v = f|v|^{p-2}v$ and $u\in \mathcal{D}_\gamma^{1,p}(\Omega)$ be a
weak supersolution of %
$$-\lap u = g|u|^{p-2}v\quad \mbox{in}\ \Omega$$
such that $\inf_{\partial \Omega}u>0$, where $g\in
L^{N_\gamma/p}(\Omega)$ such that $f\leq g$ in $\Omega$. Moreover,
assume that
\begin{equation}\label{e:conduv}
\limsup_{R\to \infty} \frac{1}{R} \int_{B_{2R}\setminus B_R} v^p
|\grad \log u|^{p-1}=0.
\end{equation}
 If $u\geq v$ on $\partial \Omega$, then
$$
u\geq v \quad in\, \Omega.
$$
\end{thm}
\begin{proof}
We prove the result for the case $p\geq2$. The case $1<p<2$ can be
proved similarly. Let $R> 2{\rm diam}(\Omega^C)$ and $t>1$. Let
$\eta_R\in C_0^\infty (B_{2R})$ be a cut-off function such that
$0\leq \eta_R \leq 1$, $\eta\equiv 1$ on $B_R$ and $|\grad
\eta_R|\leq 2/R$. Let us consider the test functions
$$
\varphi_1 = \eta_R v^{1-p} \min ((v^p-u^p)^+,t)\quad {\rm and}\ \
\varphi_2= \eta_R u^{1-p}\min ((v^p-u^p)^+, t).
$$
Testing the equations for $v$ and $u$ with $\varphi_1$ and
$\varphi_2$ respectively, we get
\begin{equation}\label{e:testfg}
\begin{split}
\langle -\lap v, \varphi_1 \rangle  - \langle -\lap u, \varphi_2
\rangle &\leq \int_\Omega f |v|^{p-2}v \varphi_1 - \int_\Omega g
|u|^{p-2}u\varphi_2\\
&= \int_{\{v^p-u^p\geq t\}} t (f-g)\eta_R - \int_{\{0\leq
v^p-u^p\leq t\}} (f-g)(v^p-u^p)\eta_R\\
& \leq 0,
\end{split}
\end{equation}
since $f\leq g$. On the other hand, by explicitly computing the
left hand side of \eqref{e:testfg} we have
\begin{equation}\label{e:lapuv}
\begin{split}
\langle -\lap v, \varphi_1 \rangle  &- \langle -\lap u, \varphi_2
\rangle\\
& =\int_{\{0\leq v^p-u^p \leq t\}} \left ( \eta_R |\grad
v|^{p-2}\grad v \cdot \grad \left ( v-\frac{u^p}{v^p}v \right ) +
\eta_R |\grad u|^{p-2}\grad u \cdot \grad \left (
u-\frac{v^p}{u^p}u \right ) \right )\\
&+\int_{\{0\leq v^p-u^p \leq t\}} \left ( \left (
v-\frac{u^p}{v^p}v \right ) |\grad v|^{p-2} \grad v \cdot \grad
\eta_R + \left ( u-\frac{v^p}{u^p}u \right ) |\grad u|^{p-2}\grad
u\cdot \grad \eta_R \right )\\
& + \int_{\{v^p-u^p \geq t\}} |\grad v|^{p-2}\grad v \cdot \grad
(t \eta_R v^{1-p}) + \int_{\{v^p-u^p\geq t \}} |\grad
u|^{p-2}\grad u \cdot \grad (-t \eta_R u^{1-p} ) \\
& =: I_1+ I_2+ I_3+I_4.
\end{split}
\end{equation}
By lemma \ref{lemma:techn}, we get that
\begin{equation}\label{e:I1}
\begin{split}
I_1 &\geq C_p \int_{\{0\leq v^p-u^p \leq t\}} \eta_R (v^p+u^p)
|\grad \log v -\grad \log u|^p\\
&\geq C_p \int_{\{0\leq v^p-u^p \leq t\}\cap B_R} (v^p+u^p) |\grad
\log v -\grad \log u|^p,
\end{split}
\end{equation}
where $C_p>0$ does not depend on $t, R$.

Now we claim that
\begin{equation} \label{e:I234}
\lim_{t,R\to \infty} I_k=0\quad {\rm for}\ k=2,3,4.
\end{equation}
Indeed, concerning $I_2$ we have
\begin{equation}
\begin{split}
|I_2|&\leq \int_{\{0\leq v^p-u^p\}} \left ( |\grad v|^{p-1}|v| +
|\grad u|^{p-1} u \right ) |\grad \eta_R| + \int_{\{0\leq
v^p-u^p\}} v^p |\grad \log u|^{p-1} |\grad \eta_R|\\
& \leq \frac{2}{R} \int_{B_{2R}\setminus B_R} \left ( |\grad
v|^{p-1}|v| + |\grad u|^{p-1}u \right )+ \frac{2}{R}
\int_{B_{2R}\setminus B_R} v^p |\grad \log u|^{p-1}=: J_1+J_2.
\end{split}
\end{equation}
By H\"{o}lder's inequality we obtain for $J_1$
\begin{equation}
\begin{split}
J_1&\leq C \left ( \int_{B_{2R}\setminus B_R} |\grad v|^{p} \right
)^{\frac{p-1}{p}} \left ( \int_{B_{2R}\setminus B_R}
|v|^{p^*_\gamma} \right )^{1/p^*_\gamma} \\
&\ \  + C \left ( \int_{B_{2R}\setminus B_R} |\grad u|^{p} \right
)^{\frac{p-1}{p}} \left ( \int_{B_{2R}\setminus B_R}
|u|^{p^*_\gamma} \right
)^{1/p^*_\gamma}\\
& =o(1)\quad {\rm as}\ R\to \infty,
\end{split}
\end{equation}
and by assumption \eqref{e:conduv} we have that
$$
J_2\to 0 \quad {\rm as}\ R\to \infty.
$$
Hence, observing also that $J_1, J_2$ are independent on $t$, we
obtain
$$
\lim_{t,R\to \infty} I_2=0.
$$
Concerning $I_3$, we observe that
\begin{equation}
I_3= t \int_{\{v^p-u^p\geq t\}} \left ( |\grad v|^{p-2} \grad v
\cdot \grad \eta_R\, v^{1-p} + (p-1) \eta_R\, v^{-p} |\grad v|^p
\right ).
\end{equation}
So
\begin{equation}
\begin{split}
|I_3|&\leq \int_{\{v^p \geq t\}} \left ( |\grad
v|^{p-1}|\grad\eta_R| v + (p-1) \eta_R |\grad v|^p \right )\\
& \leq J_1 + (p-1) \int_{\{v^p\geq t\}}|\grad v|^p,
\end{split}
\end{equation}
from which we deduce that
$$
\lim_{t,R\to \infty} I_3=0.
$$
Finally, for the term $I_4$ we observe that
\begin{equation}
\begin{split}
I_4&= \int_{\{v^p-u^p \geq t\}} \left ( t(p-1) \eta_R u^{-p}|\grad
u|^p -t u^{1-p} |\grad u|^{p-2}\grad u\cdot \grad \eta_R \right
)\\
& \geq -t \int_{\{v^p-u^p \geq t\}} |\grad \log u|^{p-1}|\grad
\eta_R|\\
& - \int_{\{v^p\geq t\}}|\grad \log u|^{p-1}|\grad \eta_R|v^p\\
&\geq -J_2,
\end{split}
\end{equation}
which tends to 0 as $t,R\to \infty$ by \eqref{e:conduv}.

So, by \eqref{e:lapuv}, \eqref{e:I1} and \eqref{e:I234}, we get
\begin{equation*}
\begin{split}
0 & \geq \limsup_{R,t\to \infty} \left ( \langle -\lap v,
\varphi_1
\rangle - \langle -\lap u, \varphi_2 \rangle \right )\\
& \geq C \int_{\{v\geq u\}} (v^p+u^p) |\grad \log v -\grad \log
u|^p.
\end{split}
\end{equation*}
Therefore
\begin{equation}
\int_{\{v\geq u\}} (v^p+u^p) |\grad \log v -\grad \log u|^p =0,
\end{equation}
from which
\begin{equation*}
\log v=\log u +C\quad {\rm on}\{z\in \Omega: v(z)\geq u(z)\},
\end{equation*}
that is
$$
v=Cu\quad {\rm on}\{z\in \Omega: v(z)\geq u(z)\},
$$
for some positive constant $C>0$. Then, by the assumption
$\inf_{\Omega} u>0$ and $v\leq u$ on $\partial \Omega$, it follows
that $C=1$. This implies
$$
v\leq u \quad {\rm in}\ \Omega.
$$
This completes the proof in the case $p\geq 2$. The case $1<p<2$
follows by analogous proof, but using inequality $ii)$ instead of
$i)$ in Lemma \ref{lemma:techn}.
\end{proof}
From the above theorem, we are able to get the lower bound
estimate for nonnegative weak solutions of our equation
\eqref{e:entire}, as nonnegative solutions of $-\lap u \geq 0$.

\medskip

\noindent \textbf{Proof of Theorem \ref{theo:asympt} - estimate
from below}. Let $u\in\mathcal{D}^{1,p}_\gamma(\RR^N)$ be a
nonnegative weak solution of \eqref{e:entire}. Then, in
particular, $u$ satisfies
\begin{equation}
-\lap u\geq 0\quad {\rm in}\ \RR^N. \label{e:superh}
\end{equation}
Hence, $u>0$  by the weak Harnack inequality for solutions to
\eqref{e:superh} (see e.g. \cite{DFZ}).

Let us compare $u$ with the function
$$
v(z)=\frac{c_1}{d(z)^{\frac{N_\gamma-p}{p-1}}}
$$
in $B_1(0)^{C}$, where $c_1:=\inf_{\partial B_1(0)}u$. As recalled
in Section \ref{sec:prelim}, a direct calculation shows that
$$
\lap v=0 \quad {\rm in}\ \RR^{N}\setminus \{0\}.
$$
Moreover, the condition \eqref{e:conduv} required in Theorem
\ref{theo:comparison} is satisfied. Indeed, we first observe that,
by adapting the arguments in Lindqvist \cite[Lemma 2.14]{Lindq2},
being $u>0$ a weak solution to \eqref{e:superh}, we have
\begin{equation}
\int_{B_{2R}\setminus B_R} |\grad \log u|^{p}\leq C
R^{N_\gamma-p}, \quad R \ {\rm large}.
\end{equation}
Then,
\begin{equation}
\begin{split}
\frac{1}{R} \int_{B_{2R}\setminus B_R} v^p |\grad \log u|^{p-1} &
\leq C R^{-1-p\frac{N_\gamma-p}{p-1} +\frac{N_\gamma}{p}} \left
(\int_{B_{2R}\setminus B_R} |\grad \log u|^p \right
)^{\frac{p-1}{p}}\\
&\leq C R^{-1-p\frac{N_\gamma-p}{p-1}+ \frac{p-1}{p}(N_\gamma-p)
+\frac{N_\gamma}{p}}\to 0 \quad \mbox{as}\ R\to \infty.
\end{split}
\end{equation}
Then, we can apply Theorem \ref{theo:comparison}, from which
$$
u(z)\geq v(z)=\frac{c_1}{d(z)^{\frac{N_\gamma-p}{p-1}}}\quad
\mbox{in}\ \ B_1(0)^C.
$$
The proof is therefore complete.

\section{Existence results for Brezis-Nirenberg type problems}
\label{sec:BNexist}

\medskip

This section is devoted to the proof of Theorem
\ref{boundedmainresult1} and Theorem \ref{boundedmainresult2}. To
this aim we look for critical points of the energy functional
$J_\lambda: \mathring{W}^{1,p}_{\gamma}(\Omega) \rightarrow
\mathbb{R}$ defined by

$$
J_\lambda(u):=\frac{1}{p} \int_{\Omega} \left|\nabla_{\gamma} u
\right|^p-\frac{\lambda}{q} \int_{\Omega} \abs{u}^q  d
z-\frac{1}{p_\gamma^*} \int_{\Omega}\abs{u}^{p_\gamma^*} d z,
$$
where $p\leq q< p_\gamma^*$. The existence results, as in the
ordinary $p$-Laplacian case, will follow by an application of the
Mountain pass Theorem, exploiting the compactness properties of
$J_\lambda$ under a suitable threshold.

We begin by stating that $J_\lambda$ satisfies the Palais-Smale
condition under a certain level. For the Euclidean $p$-Laplacian
case, see \cite[Lemma 2.3]{GP2}, see also \cite[Lemma 1]{AG}; a
recent fractional $p$-Laplacian version can be found in
\cite{PSY}. In the Grushin context, it appears in \cite[Lemma
6.1]{MMS} for the case $q=p$.

\begin{lem}\label{lemma:thresh}
    Let $\{u_k\}\subset \mathring{W}^{1,p}_{\gamma}(\Omega)$ be a
    sequence such that
  \begin{equation*}
    J_\lambda(u_k)\to c\in \left (0,
    \frac{\mathcal{S}_{\gamma,p}^{\frac{N_\gamma}{p}}}{N_\gamma}\right )
\end{equation*}
    and
\begin{equation*}
    J_\lambda'(u_k)\to 0 \ \ \mbox{in}\, \left ( \mathring{W}^{1,p}_{\gamma} (\Omega) \right
    )^{*}.
\end{equation*}
    Then, there exists $u\in \mathring{W}^{1,p}_{\gamma}(\Omega)\setminus \{0\}$ such
    that $u_k\rightharpoonup u$, up to a subsequence, and
    $J_\lambda'(u)=0$.
\end{lem}
\begin{proof}
 We give here a sketch of the proof, following the
outline in \cite{AG}. Let $f(s)=\lambda|s|^{q-2}s
+|s|^{p^*_\gamma-2}s$ and $F(s)=\int_0^s f(t)\,dt$. Since it holds
that
$$
\exists \,\theta\in \left (0, \frac{1}{p} \right ), \exists\,
\bar{s}>0\,\, \mbox{such that}\, F(s)\leq \theta f(s)s, \ \forall
|s|\geq \bar{s},
$$
then $\{u_k\}$ is bounded (see e.g. \cite{R}) and there exists
$u\in \mathring{W}^{1,p}_{\gamma}(\Omega)$ such that, up to a
subsequence, $u_k\rightharpoonup u$. Moreover, $J'_\lambda (u)=0$
by the weak continuity of $J'_\lambda$.

Assume by contradiction that $u=0$. Since the term $\|u_k\|_q^{q}$
is subcritical, from $J'_\lambda(u_k)[u_k]=o(1)$ we infer that
\begin{equation} \label{e:step}
\|u_k\|^p -\|u_k\|_{p_\gamma^*}^{p_\gamma^*} =o(1).
\end{equation}
Then, by the Sobolev inequality \eqref{e:Sobol}, we get
$$
o(1)\geq \|u_k\|^p\left ( 1-S_{\gamma,p}^{-\frac{p^*_\gamma}{p}}
\|u_k\|^{p^*_\gamma -p} \right ).
$$
Now, observe that $\|u_k\|\nrightarrow 0$ since $c>0$. Hence, the
last inequality implies
$$
\|u_k\|^p \geq S_{\gamma,p}^{N_\gamma/p} + o(1),
$$
and by \eqref{e:step} we deduce
\begin{equation}
\begin{split}
J_\lambda (u_k) &= \frac{1}{N_\gamma}\|u_k\|^p+
\frac{1}{p_\gamma^*} \left
(\|u_k\|^p-\|u_k\|_{p_\gamma^*}^{p^*_\gamma} \right ) + o(1)\\
& \geq \frac{1}{N_\gamma} S_{\gamma,p}+ o(1),
\end{split}
\end{equation}
which contradicts the assumption $c<\frac{1}{N_\gamma}
S^{N_\gamma/p}$.
\end{proof}

It is easily seen that, under the assumptions of Theorem
\ref{boundedmainresult1} and \ref{boundedmainresult2}, the
functional $J_\lambda$ satisfies the Mountain pass geometry
requirements.

\begin{lem}\label{MP-geometry0}
    The functional $J_{\lambda}$, for $q=p$ and $0<\lambda<\lambda_1$, or $p<q<p^*_\gamma$ and
$\lambda>0$, satisfies
    \begin{enumerate}
        \item[(i)] There exists positive constants $\delta$ and $R$ such that $J_{\lambda}(v) \geq \delta $ for all
        $v\in \mathring{W}^{1,p}_{\gamma}(\Omega)$, $\norm{v} =
        R$;
        \item[(ii)] There exists $v_{1} \in \mathring{W}^{1,p}_{\gamma}(\Omega)$ with $\norm{v_{1}} >R$ and $J_{\lambda}(v_{1})<0.$
    \end{enumerate}
\end{lem}
\begin{proof}
    (i) If $q=p$ and $0<\lambda<\lambda_1$, by the definition of
    $\lambda_1$ and the Sobolev embedding, we get
    \begin{align*}
        J_{\lambda}(u)
        \geq \frac{1}{p} \left ( 1-\frac{\lambda}{\lambda_1} \right ) \|u\|^p -
        \frac{C}{p_{\gamma}^{*}} \norm{u}^{p_{\gamma}^{*}}.
    \end{align*}
   So, $0$ is a strict local minimum and (i) follows. Analogously,
    if $p<q<p^*_\gamma$ and
$\lambda>0$, the continuous embedding $
\mathring{W}^{1,p}_{\gamma}(\Omega) \hookrightarrow L^{q}(\Omega)
$ for $ q\leq p_{\gamma}^{*}$ implies that
    \begin{align*}
        J_{\lambda}(u)& =\frac{1}{p} \int_{\Omega} \left|\nabla_{\gamma} u\right|^p-\frac{\lambda}{q} \int_{\Omega} \abs{u}^q  d z-\frac{1}{p_\gamma^*} \int_{\Omega}\abs{u}^{p_\gamma^*} d z \\
        & \geq \frac{1}{p} \norm{u}^p-\lambda \frac{C_1}{q} \norm{u}^q-\frac{C_2}{p_{\gamma}^{*}} \norm{u}^{p_{\gamma}^{*}}.
    \end{align*}
    form which (i) easily follows.

(ii) The condition obviously holds since the leading term at
$\infty$ is the critical term.


\end{proof}

 Now, in view of Lemma \ref{MP-geometry0}, the \textit{mountain-pass}
level defined by
\begin{align}\label{vcritical0}
    c_\lambda= \inf_{\alpha \in \mathcal{A}} \max_{t \in[0,1]} J_{\lambda}(\alpha (t)),
\end{align}
where%
$$\mathcal{A}= \left\{\alpha \in C\left([0,1],
\mathring{W}^{1,p}_{\gamma}\left(\Omega \right) \right) \mid
\alpha(0)=0, J(\alpha(1))<0 \right\}
$$
is well-defined, and by Theorem 2.2 in \cite{Bre} (see also
\cite[Lemma 3.1]{GP2}), $J_\lambda$ admits a PS-sequence at level
$c_\lambda$; moreover, such a sequence may be chosen in the cone
of nonnegative functions since $J_\lambda(|u|)\leq J_\lambda(u)$,
for all $u\in \mathring{W}^{1,p}_{\gamma} (\Omega).$ We will find
our solutions by proving that, under the assumptions (on $\lambda$
and $q$), respectively of Theorem \ref{boundedmainresult1} and
\ref{boundedmainresult2}, $c_\lambda$ stays below the compactness
threshold  given by Lemma \ref{lemma:thresh}.

To this aim, following the usual Brezis-Nirenberg argument, a
crucial role is played by suitable truncations (out of $\Omega$)
of the $p$-Sobolev extremals. Let's introduce this family, arguing
as in \cite{Loiu00}.

Let $U>0$ be a fixed extremal function for inequality
 \eqref{e:Sobol1}. We can assume, up to a normalization, that $\|\grad U\|_p^p=\|U\|_{p_\gamma^*}^{p_\gamma^*}=S_{\gamma, p}^{N_\gamma/p}$.
 For $\epsilon>0$, define the rescalings
 \begin{equation}\label{e:Ueps}
U_\epsilon (z)= \epsilon^{-(N_\gamma-p)/p}
U(\delta_{1/\epsilon}(z)), \quad z \in \RR^N.
 \end{equation}
Of course, $U_\epsilon$ are also minimizers and verify
\begin{equation}\label{e:normal}
\|\grad
U_\epsilon\|_p^p=\|U_\epsilon\|_{p_\gamma^*}^{p_\gamma^*}=S_{\gamma,p}^{N_\gamma/p}.
\end{equation}
Now, due to $\Omega \cap \{x=0\}\neq \emptyset$, it is not
restrictive to suppose $0\in \Omega$, due to the translation
invariance with respect to the $y$-variable. Let $R>0$ be such
that $B_R\Subset \Omega$ and let $\varphi \in C_0^\infty (B_R)$ be
a cut-off function, $0\leq \varphi \leq 1$, $\varphi \equiv 1$ in
$B_{R/2}$. Define
\begin{equation}
u_{\epsilon}(z) := \varphi(z)U_{\epsilon}(z), \ \ \epsilon
>0.\label{ueps}
\end{equation}
Reasoning as in \cite{Loiu00}, by means of the asymptotic
estimates proved in Theorem \ref{theo:asympt}, we are able to
prove the following asymptotic expansions:
\begin{lem} \label{lemma:expans}
The following estimates hold for $\epsilon \to 0$
\begin{align}\label{e:expans}
 \|\grad u_\epsilon \|_p^p &\leq
S_{\gamma,p}^{{N_\gamma}/p} + C \epsilon^{(N_\gamma-p)/(p-1)},\\
\quad \|u_\epsilon\|_{p_\gamma^*}^{p_\gamma^*}&\geq
S_{\gamma,p}^{{N_\gamma}/p}-C \epsilon^{N_\gamma/(p-1)}\\
 \norm{u_{\epsilon}}_{p}^{p}  &= \begin{cases}
    C \epsilon^{p} + O\left(\epsilon^{\frac{N_{\gamma}-p}{p-1}}\right) \text{ if } N_{\gamma} > p^2, \\
    C\epsilon^{p}\abs{\ln \epsilon} + O\left(\epsilon^p\right) \text{ if } N_{\gamma} = p^2, \\
    C\epsilon^{\frac{N_{\gamma}-p}{p-1}} + O(\epsilon^p)\ \ \text{ if } N_{\gamma} < p^2.
 \end{cases}
\end{align}
\end{lem}
\begin{proof}
Let us begin with the first estimate. First observe that
\begin{equation}
\begin{split} \label{e:grad1}
    \int_{\Omega} \abs{\nabla_{\gamma}u_{\epsilon}}^p dz&=
    \int_{\Omega} |\grad (\varphi U_\epsilon)|^p dz \\
& = \int_\Omega \left | \varphi \grad U_\epsilon+
U_\epsilon \grad \varphi  \right |^{p} dz\\
& \leq \left ( \|\varphi \grad U_\epsilon \|_p + \|U_\epsilon
\grad \varphi\|_p \right )^p\\
& \leq \left ( \|\grad U_\epsilon\|_p + \|U_\epsilon \grad
\varphi\|_p \right )^p.
\end{split}
    \end{equation}
For the last term in the above estimate we have
    \begin{equation}\label{e:grad2}
    \begin{split}
\int_{\Omega} U_\epsilon^p |\nabla_\gamma \varphi|^p  \,dz& \leq
C\,\int_{{B_R}\setminus {B_{R/2}}} U_\epsilon^p \, dz
=C\,\epsilon^{p-N_\gamma}\int_{{B_R}\setminus {B_{R/2}}}
 \, U^p(\delta_{\frac{1}{\epsilon}}z)\, dz\,\\
&= C \epsilon^p \int_{B_{R/\epsilon}\setminus B_{R/2\epsilon}} U^p\,dz \\
&\leq C \epsilon^p \int_{R/2\epsilon<d(z) < R/\epsilon}
\frac{1}{\,\,\,d(z)^{\frac{p(N_\gamma-p)}{p-1}}}\,dz\\
&= C \epsilon^p \int_{R/2\epsilon}^{R/\epsilon}
\frac{\rho^{N_\gamma-1}}{\,\,\,\rho^{\frac{p(N_\gamma-p)}{p-1}}}\,d\rho\\
&= O(\epsilon^{(N_\gamma-p)/(p-1)}),
\end{split}
\end{equation}
where we have used the estimate from above on $U$ for $d(z)$
large. Hence, from \eqref{e:grad1}, \eqref{e:grad2} and
\eqref{e:normal}, estimate \eqref{e:expans} follows. Next, we
compute
\begin{equation}
\begin{split}
 \int_{\Omega} u_{\epsilon}^{p_{\gamma}^{*}} dz &=   \int_{\Omega} \varphi^{p_{\gamma}^{*}} U_{\epsilon}^{p_{\gamma}^{*}} dz
    = \int_{\Omega}  U_{\epsilon}^{p_{\gamma}^{*}} dz + \int_{\Omega} \left( \varphi^{p_{\gamma}^{*}}-1\right) U_{\epsilon}^{p_{\gamma}^{*}} dz \\
    & = \int_{\mathbb{R}^N}  U_{\epsilon}^{p_{\gamma}^{*}} dz -\int_{\mathbb{R}^N \setminus \Omega}  U_{\epsilon}^{p_{\gamma}^{*}} dz
    -\int_{\Omega} \left( 1-\varphi^{p}\right) U_{\epsilon}^{p_{\gamma}^{*}} dz \\
    &= S_{\gamma,p}^{{N_\gamma}/p}+ \alpha(\varphi, \epsilon),
\end{split}
\end{equation}
where $\alpha(\varphi,\epsilon)= O(\epsilon^{N_\gamma/(p-1)})$.
Indeed,
\begin{align}
0 & \leq \int_\Omega (1-\varphi^p) U_\epsilon^{p_\gamma^*}
\,dz\leq \int_{d(z)>R/2} U_\epsilon^{p_\gamma^*}\, dz =
\int_{d(z)>
R/2\epsilon} U^{p_\gamma^*}\, dz \nonumber \\
{} & \leq C\,\int_{d(z)> R/2\epsilon} \frac{1}{d(z)^{\frac{p
N_\gamma}{p-1}}}\,dz = C\,\int_{R/2\epsilon}^{\infty}
\frac{1}{\rho^{N_\gamma+p-1}}\,{\rm d}\rho \nonumber \\
{}& = O(\epsilon^{N_\gamma/(p-1)}), \nonumber
\end{align}
and an analogous estimate holds for the other term in
$\alpha(\varphi,\epsilon)$.

Finally, we estimate
\begin{align*}
\int_{\Omega} u_{\epsilon}^{p} dz & \geq \int_{B_{\frac{R}{2}}}
U_{\epsilon}^{p} dz  = \int_{B_{\frac{R}{2}}}
\epsilon^{p-N_{\gamma}} U^{p} \left(\delta_{1/\epsilon}(z) \right)
dz  = \epsilon^{p} \int_{B_{\frac{R}{2\epsilon}}} U^{p}(z) dz
\\ &= \epsilon^{p} \left( \int_{B_{R_0}} U^{p}(z) dz +
\int_{B_{\frac{R}{2\epsilon}} \setminus B_{R_0}} U^{p}(z) dz
\right) \\ & \geq C \epsilon^{p} \left( 1 +
\int_{B_{\frac{R}{2\epsilon}} \setminus B_{R_0}}
\frac{1}{d(z)^{\frac{p(N_{\gamma}-p)}{p-1}}} \right) dz \\ & \geq
C \epsilon^{p} \left( 1 + \int_{R_0}^{R/2\epsilon}
\frac{\rho^{N_{\gamma}-1}}{\rho^{\frac{p(N_{\gamma}-p)}{p-1}}} d
\rho \right) \\ 
&= \begin{cases}
    C \epsilon^{p} + O\left(\epsilon^{\frac{N_{\gamma}-p}{p-1}}\right) \text{ if } N_{\gamma} > p^2, \\
    C\epsilon^{p}\abs{\ln \epsilon} + O\left(\epsilon^p\right) \text{ if } N_{\gamma} = p^2,
    \\ C\epsilon^{\frac{N_{\gamma}-p}{p-1}} + O(\epsilon^p) \text{ if } N_{\gamma} <
    p^2,
\end{cases}
\end{align*}
where we have used the estimate from below on $U$.
\end{proof}

\noindent \textbf{Proof of Theorem \ref{boundedmainresult1}:} (i)
Let $N\geq p^2$ and $0<\lambda<\lambda_1$. Since
$$
J_\lambda (tu_\epsilon)\to -\infty \quad \mbox{as}\ t\to \infty,
$$
fix $R_0>0$ so large that $J_\lambda (R_0 u_\epsilon )<0$. Then,
\begin{equation}
\begin{split}
c_\lambda\leq \max_{t\in[0,1]} J_\lambda (tR_0u_\epsilon)
&=\frac{1}{N_\gamma} \left ( \frac{\|\grad u_\epsilon \|_p^p
-\lambda \|u_\epsilon\|_p^p}{\|u_\epsilon\|_{p^*_\gamma}^p} \right
)^{N_\gamma/p} =: \frac{1}{N_\gamma}Q_\lambda
(u_\epsilon)^{N_\gamma/p}.
\end{split}
\end{equation}
Now, by Lemma \ref{lemma:expans}
\begin{align*}
    Q_{\lambda} (u_{\epsilon}) & \leq \begin{cases}
 \mathcal{S}_{\gamma,p}+ C\epsilon^{(N_\gamma-p^2)/(p-1)} -\lambda
\epsilon^{p},& N_\gamma>p^2\\
\mathcal{S}_{\gamma,p}+ C\epsilon^p -\lambda C \epsilon^{p}
\abs{\ln \epsilon}, & N_\gamma=p^2.
\end{cases}
\end{align*}
Hence, $Q_{\lambda} (u_{\epsilon})< S_{\gamma,p}$ if $\epsilon>0$
is sufficiently small and so
$$
c_\lambda< \frac{1}{N_\gamma}S_{\gamma,p}^{N_\gamma/p},
$$
hence we can conclude.\\
 (ii)   For
$N_\gamma<p^2$, let us consider an eigenfunction $\varphi_1$
associated with $\lambda_1=\lambda_1(\Omega)$. Then
$Q_{\lambda_1}(\varphi_1)=0$; by continuity, there exists
$\lambda_*<\lambda_1$ such that
$Q_\lambda(\varphi_1)<\mathcal{S}_{\gamma,p}$ for
$\lambda_*<\lambda<\lambda_1.$\\

We proceed to prove Theorem \ref{boundedmainresult2}.

\begin{prop}\label{vcriticalbound0}
If $\max\{p, p^*_\gamma (1-\frac{1}{p}), p^*_\gamma-
\frac{p}{p-1}\}< q<p^*_\gamma$, the number $c_{\lambda}$ defined
by
    (\ref{vcritical0}) satisfies
    \begin{align}\label{condclambda}
        c_{\lambda} < \frac{1}{N_\gamma}
        \mathcal{S}^{\frac{N_\gamma}{p}}_{\gamma},\quad
        \forall \lambda>0.
    \end{align}
\end{prop}
\begin{proof}
    Let us consider
    \begin{equation} \label{wA}
        w_{\epsilon}:= \frac{u_{\epsilon}}{\norm{u_{\epsilon}}}_{L^{p_{\gamma}^{*}}(\Omega)},
    \end{equation}
    where the functions $u_{\epsilon}$ are defined as in \eqref{ueps}.
    Since $J_{\lambda}(tw_{\epsilon}) \to -\infty$ as $t \to +\infty,$ there exists $R>0$ such that $J_{\lambda}(Rw_{\epsilon})<0.$ Now, defining $v_{1}= Rw_{\epsilon},$ and using Lemma \ref{MP-geometry0}, we can write
    $$
    c_{\lambda}=\inf_{\alpha \in \mathcal{A}} \max_{t \in[0,1]} J_{\lambda}(\alpha (t)) \leq \max_{t >0}J_{\lambda}(tw_{\epsilon}).
    $$
    Moreover, reasoning as in \cite[Lemma 5.3]{AGLT}, it is easy to see
    that there exists a unique $t_{\epsilon}>0$ such that
    $$
    J_{\lambda}(t_{\epsilon}w_{\epsilon})=\max_{t > 0} J_{\lambda}(tw_{\epsilon}).
    $$
    So, in order to prove the lemma, it is sufficient to show that
    $$
    J_{\lambda}(t_{\epsilon}w_{\epsilon}) < \frac{1}{N_\gamma} \mathcal{S}^{\frac{N_\gamma}{p}}_{\gamma,p}.
    $$
Since $\norm{w_{\epsilon}}_{L^{p_{\gamma}^{*}}(\Omega)}=1,$ we
have for $t>0$
    \begin{align*}
        J_{\lambda}\left( tw_{\epsilon} \right) &=  \frac{t^p}{p} \int_{\Omega} \abs{\nabla_{\gamma}w_{\epsilon}}^p dz -
        \lambda \frac{t^q}{q} \int_{\Omega}  \abs{w_\epsilon}^q dz- \frac{t^{p_{\gamma}^{*}} }{{p_{\gamma}^{*}}} \int_{\Omega}  \abs{w_{\epsilon}}^{p_{\gamma}^{*}} dz\\
        &= \left(\frac{t^p}{p} \int_{\Omega} \abs{\nabla_{\gamma}w_{\epsilon}}^2dz-\frac{t^{p_{\gamma}^{*}} }{{p_{\gamma}^{*}}} \right) -
        \lambda \frac{t^q}{q} \int_{\Omega}  \abs{w_{\epsilon}}^q dz.\\
    \end{align*}
    Define $$h(t):= \frac{t^p}{p} \int_{\Omega} \abs{\nabla_{\gamma}w_{\epsilon}}^pdz-\frac{t^{p_{\gamma}^{*}} }{{p_{\gamma}^{*}}}.$$
    It is easy to check that $h$ achieves its  maximum at point
    $s_{\epsilon}:= \left( \int_{\Omega} \abs{\nabla_{\gamma}w_{\epsilon}}^p dz
    \right)^{\frac{1}{p_{\gamma}^{*}-p}}$ and it is
    increasing on $[0, s_\epsilon]$. Moreover, from
    $$
    0=J'_\lambda(t_\epsilon)=t_{\epsilon}^{p-1} \left
    (\int_{\Omega}\abs{\nabla_{\gamma}w_{\epsilon}}^p \,dz-t_\epsilon^{p^*_\gamma-p}-\lambda
    t_\epsilon^{q-p}\int_\Omega |w_\epsilon|^q dz \right ),
    $$
    we have
    $$
    \int_{\Omega} \abs{\nabla_{\gamma}w_{\epsilon}}^p \,dz = t_\epsilon^{p^*_\gamma-p}+\lambda
    t_\epsilon^{q-2}\int_\Omega |w_\epsilon|^q dz \geq
    t_\epsilon^{p_\gamma^*-p}.
    $$
    Hence $$t_{\epsilon} \leq \left(\int_{\Omega} \abs{\nabla_{\gamma} w_{\epsilon}}^p dz \right)^{\frac{1}{p_{\gamma}^{*}-p}}=s_\epsilon.$$
    Therefore, also taking into account that, by
    Lemma \ref{lemma:expans}, it holds
    $$
    \int_{\Omega} \abs{\nabla_{\gamma} w_{\epsilon}}^p dz= \mathcal{S}_{\gamma,p} +
    O(\epsilon^{\frac{N_\gamma-p}{p-1}}),
    $$
    we get that
    \begin{align}\label{e:condq}
        J_{\lambda}(t_{\epsilon}w_{\epsilon}) &=h(t_{\epsilon})-\frac{\lambda}{q}t_\epsilon^q \int_\Omega |w_\epsilon|^q\, dz \nonumber \\
        &\leq h(s_\epsilon)-\frac{\lambda}{q}t_\epsilon^q \int_\Omega |w_\epsilon|^q\, dz \nonumber \\
        &=\frac{1}{N_\gamma}\left (\int_{\Omega}
        \abs{\nabla_{\gamma}w_{\epsilon}}^p\,dz \right
        )^{\frac{N_\gamma}{p}}-\frac{\lambda}{q}t_\epsilon^q \int_\Omega
        |w_\epsilon|^q\,
        dz \nonumber \\
        &\leq\frac{1}{N_{\gamma}}\mathcal{S}_{\gamma}^{\frac{N_{\gamma}}{p}}+
        O(\epsilon^{\frac{N_\gamma-p}{p-1}})  -
        \lambda C_{\epsilon} \int_{\Omega} |w_{\epsilon}|^q dz,
    \end{align}
    where $C_{\epsilon}=\frac{t_{\epsilon}^q}{q}.$ We observe that there is a positive constant
    $C_{0}$ such that $$C_{\epsilon} \geq C_{0}>0, \forall
    \epsilon>0.$$
    Otherwise, we could find a sequence $\epsilon_k \rightarrow 0,$ as $k \rightarrow \infty$ such
    that $t_{\epsilon_k} \rightarrow 0$ as $k \rightarrow \infty.$ So, up to a subsequence, we should have $t_{\epsilon_k}w_{\epsilon_k}
\rightarrow 0,$ as $k \rightarrow \infty.$
    Therefore, $$0<c_{\lambda} \leq \max_{t \geq 0} J_{\lambda}(tw_{\epsilon_k}) =J_{\lambda}(t_{\epsilon_k}w_{\epsilon_k})=
    J_{\lambda}(0)=0,$$ which is a contradiction.

Moreover, reasoning as in Lemma \ref{lemma:expans}, we have the
    following estimates, as $\epsilon \to 0$:
    \begin{align}
        \int_\Omega |w_\epsilon|^q \, dz & \geq
        \begin{cases}
            O(\epsilon^{N_\gamma- \frac{(N_\gamma-p)}{p}q})  & {\rm
                if}\,\, p^*_\gamma  (1-\frac{1}{p} )<q<p_\gamma^*\\
            O(\epsilon^{N_\gamma- \frac{(N_\gamma-p)}{p}q}|\ln \epsilon|) & {\rm
                if}\,\,q=p^*_\gamma (1-\frac{1}{p}) \\
            O(\epsilon^{\frac{(N_\gamma-p)}{p}q}) & {\rm if}\,\, 1\leq
            q<p^*_\gamma (1-\frac{1}{p}), \label{e:brez4}
        \end{cases}
    \end{align}
where note that
$$
p^*_\gamma \left (1-\frac{1}{p} \right )= \frac{N_\gamma
(p-1)}{N_\gamma-p}.
$$
    So, by \eqref{e:condq} and
       \eqref{e:brez4}, and taking into account that
    $$
    N_\gamma-\frac{N_\gamma-p}{p}q<\frac{N_\gamma -p}{p-1}\ \
    \mbox{for}\ \ q>\frac{(N_\gamma p -2 N_\gamma +p)p}{(p-1)
(N_\gamma-p)}=p^*_\gamma- \frac{p}{p-1},$$ we have that, if
$\max\{p,\, p^*_\gamma (1-\frac{1}{p}),\, p^*_\gamma-
\frac{p}{p-1}\}<q<p_\gamma^*$, then
    $$J_\lambda(t_\epsilon w_\epsilon) < \frac{\mathcal{S}^{\frac{N_\gamma}{p}}_{\gamma,p}}{N_\gamma}, \text{ for } \epsilon>0 \text{ sufficiently small},$$
    for any $\lambda>0$.  This proves our aim.
\end{proof}
\noindent \textbf{ Proof of Theorem \ref{boundedmainresult2}:} (i)
Let $N_\gamma\geq p^2$. Then, $p^*_\gamma- \frac{p}{p-1}\leq
p^*_\gamma (1-\frac{1}{p})\leq p$, as observed in \cite[Remark
3.4]{GP2}. So by Lemma \ref{vcriticalbound0}, we get that
condition \eqref{condclambda} holds for any $p<q<p_\gamma^*$.
Therefore, we get a nontrivial solution to problem
\eqref{e:mainp1} for any $\lambda>0$, for the whole range of $q$. \\
(ii) Let $N_\gamma<p^2$. In this case, $p<p^*_\gamma
(1-\frac{1}{p})<p^*_\gamma- \frac{p}{p-1}$, so condition
\eqref{condclambda} in Lemma \ref{vcriticalbound0} holds for
$p^*_\gamma- \frac{p}{p-1}<q< p_\gamma^*$. Hence, for this range
of $q's$ we get a nontrivial solution to problem \eqref{e:mainp1}
for any $\lambda>0$.

For $p<q\leq p^*_\gamma- \frac{p}{p-1}$, the estimate is
insufficient to get \eqref{condclambda}. However, we can get that
there exists $\lambda^*>0$ such that
\begin{equation}\label{e:condlarge}
c_\lambda <\frac{1}{N_\gamma}
\mathcal{S}_{\gamma,p}^{\frac{N_\gamma}{p}}, \quad \forall \lambda
\geq \lambda_{*}.
\end{equation}
Indeed, for any fixed nonnegative function $u \in
\mathring{W}^{1,p}_{\gamma}\left(\Omega \right) \setminus \{0\}$,
it holds
$$
M_\lambda:=\max_{t \geq 0}\left\{
\frac{t^p}{p}\int_{\Omega}|\nabla_{\gamma} u|^p\,dz-\lambda \frac{
    t^q}{q}\int_{\Omega}|u|^q\,dz
\right\}=\frac{1}{\lambda^{\frac{p}{q-p}}}\left(\frac{1}{p}-\frac{1}{q}\right)\left(\frac{\|u\|^p}{\,\,\,\|u\|_q^{q}}
\right)^{\frac{q}{q-p}} \to 0 \quad \mbox{as} \quad \lambda \to
+\infty.
$$
Thereby, by definition of $c_\lambda$,
$$
0<c_\lambda \leq \max_{t \geq 0}J_\lambda(tu) \leq M_\lambda,
\quad \forall \lambda >0,
$$
from which \eqref{e:condlarge} follows. This concludes the proof.
\medskip



\begin{thebibliography}{9}

\bibitem{AFL} L. Abatangelo, A. Ferrero, P. Luzzini, \textit{On
solutions to a class of degenerate equations with the Grushin
operator}, preprint 2024 https://arxiv.org/abs/2410.12637

\bibitem{Alv} C. O. Alves, A. R. F. de Holanda, \textit{A Berestycki-Lions type result for a class of degenerate elliptic
problems involving the Grushin operator,} Proc. R. Soc. Edinb.,
128 (2022), DOI:10.1017/prm.2022.43.

\bibitem{Alv new} C. O. Alves, A. R. F. de Holanda,
\textit{Existence of the solution for a class of the semilinear
degenerate elliptic equation involving the Grushin operator in
$\RR^2$: The interaction between Grushin's critical exponent and
exponential growth}, Math. Nachr. (2023),
DOI:10.1002/mana.202300171

\bibitem{AGLT} C. O. Alves, S. Gandal, A. Loiudice,
J. Tyagi, \textit{A Br\'{e}zis-Nirenberg Type Problem for a Class
of Degenerate Elliptic Problems Involving the Grushin Operator},
J. Geom. Anal. (2024), 34(2), 52.

\bibitem{AG} G. Arioli, F. Gazzola, \textit{Some results on p-Laplace equations with a critical growth
term}, Differential and Integral Equations (1998), 11(2):311-326.


\bibitem{Baou} M. S. Baouendi, \textit{Sur une Classe d'Operateurs Elliptiques Degener\'{e}s,} Bull. Soc. Math. France 95 (1967), 45-87.

\bibitem{Beck} W. Beckner, \textit{On the Grushin operator and
hyperbolic symmetry}, Proc. Amer. Math. Soc. 129 (2001), no. 4,
1233-1246.

\bibitem{B} T. Bieske, \textit{Fundamental solutions to the p-Laplace equation in a class of grushin vector fields},
 Electronic Journal of Differential Equations, Vol. 2011 (2011), No. 84, pp. 1-10.

\bibitem{BG} T. Bieske, J. Gong,  \textit{The p-Laplace equation on a class of
Grushin-type spaces}, Proc. Amer . Math. Soc. \textbf{134} (2006),
no. 12, 3585-3594.

\bibitem{Squass} L. Brasco, S. Mosconi, M. Squassina,
\textit{Optimal decay of extremals for the fractional Sobolev
inequality}, Calc. Var. Partial Differential Equations 55 (2015),
55-23.


\bibitem{Bre} H. Brezis and L. Nirenberg, \textit{Positive solutions of nonlinear elliptic equations involving critical exponents,}
Comm. Pure Appl. Math. 36 (1983) 437-477.
%
\bibitem{CDG1} L. Capogna, D. Danielli, N. Garofalo, \textit{An
embedding theorem and the Harnack inequality for nonlinear
subelliptic equations}, Comm. Partial Differential Equations 18
no. 9-10 (1993), 1765-1794

\bibitem{Hua Chen} H. Chen, H.G. Chen, J.N. Li, \textit{Sobolev
inequality and its extremal functions for homogeneous
H\"{o}rmander vector fields}, preprint 2025
https://doi.org/10.48550/arXiv.2506.16125v2

\bibitem{Chen} H. Chen, X. Liao, M. Zhang, \textit{Dirichlet problem
for a class of nonlinear degenerate elliptic operators with
critical growth and logarithmic perturbation}, Calc. Var. (2024)
63:116.


\bibitem{DM} L. D'Ambrosio, E. Mitidieri, \textit{A priori estimates and reduction
principles for quasilinear elliptic problems and applications},
Adv. Differential Equations 17 (2012) 935-1000.

\bibitem{DFZ} G. Di Fazio, M. S. Fanciullo, P. Zamboni,
\textit{Harnack inequality and continuity of weak solutions for
doubly degenerate elliptic equations}, Math. Z. (2019)
292:1325-1336.

\bibitem{Dou} J. Dou, L. Sun, L. Wang, M. Zhu, \textit{Divergent operator with degeneracy and related sharp inequalities,}
J. Functional Anal. 282 (2022), no. 2, Paper No. 109294, 85 pp.

\bibitem{Duo0} A. T. Duong, N. T. Nguyen, \textit{Liouville type theorems for elliptic equations involving Grushin
operator and advection,} Electron. J. Differ. Equ. 2017 (2017),
1-11.

\bibitem{Duo1} A. T. Duong, Q. H. Phan, \textit{Liouville type theorem for nonlinear
elliptic system involving Grushin operator,} J. Math. Anal. Appl.
454 (2017), 785-801.


\bibitem{FL1} B. Franchi, E. Lanconelli, \textit{Une
metrique associ\'{e}e \`{a} une classe d'op\'{e}rateurs
elliptiques d\'{e}g\'{e}n\'{e}r\'{e}s, Proceedings of the meeting
"Linear partial and pseudodifferential operators"}, Rend. Sem.
Mat. Univ. e Politec. Torino 1982, 105-114.

\bibitem{FL2} B. Franchi, E. Lanconelli, \textit{An embedding theorem for Sobolev Spaces related to non-smooth vector fields and
 Harnack inequality}, Comm. Partial Differential Equations 9 (1984) 1237-1264.

 \bibitem{FL3} B. Franchi, E. Lanconelli, \textit{H\"{o}lder regularity theorem for a class
of linear nonuniformly elliptic operators with measurable
coefficients}, Ann. Sc. Norm. Sup Pisa Cl Sci. (4) 1983, 523-541.

\bibitem{GP1} J. Garc\`{i}a Azorero, I. Peral
Alonso, \textit{Existence and nonuniqueness for the $p$-Laplacian:
nonlinear eigenvalues}, Comm. Partial Differential Equations
\textbf{12}, 1389-1430 (1987)

\bibitem{GP2} J. Garc\`{i}a Azorero, I. Peral Alonso,
\textit{Multiplicity of solutions for elliptic problems with
critical exponent or with a nonsymmetric term}, Trans. Amer. Math.
Soc. \textbf{323}, 2, 877-895 (1991)


\bibitem{G} N. Garofalo, \textit{Unique continuation for a class
 of elliptic operators which degenerate  on a manifold of arbitrary codimension},
Journal of Differential Equations
 104(1993) 117-146.


\bibitem{Graf} L. Grafakos, Classical Fourier Analysis.
2nd Edition, Graduate Text in Mathematics 249, Springer, New York,
2008.

\bibitem{Gr1} V. V. Grushin, \textit{On a class of hypoelliptic operators,} 1970 Math. USSR Sb. 12 458.

\bibitem{Gr2} V. V. Grushin, \textit{On a class of hypoelliptic pseudodifferential operators degenerate on a submanifold,} 1971 Math. USSR Sb. 13 155.

\bibitem{GJ} Y. Guo, Y. Jiang, \textit{Harnack inequality for
subelliptic $p$-Laplacian equations of Schr\"{o}dinger type},
Journal of Inequalities and Applications 2013, 2013:160.

\bibitem{GL} C. Gutierrez, E. Lanconelli, \textit{Maximum
principle, nonhomogeneous Harnack inequality, and Liouville
theorems for $X$-elliptic operators}, Comm. Partial Differential
Equations 28 (11-12) (2003), 1833-1862.

\bibitem{HY} J. Z. Huang, X. X. Yang, \textit{$P$-Laplace equation for
the Grushin type operator}, Acta Math. Sinica, 2023, No. 5, pp.
923-938.

\bibitem{Kog} A. E. Kogoj, E. Lanconelli, \textit{On semilinear $\Delta_{\lambda}$-Laplace equation,} Nonlinear Anal. 75 (2012), 4637-4649.

\bibitem{KLV} A. E. Kogoj, M. E. Lima, A. Viana, \textit{On the heat equation
involving a Grushin operator in Marcinkiewicz spaces}, Rev Mat
Complut (2025). https://doi.org/10.1007/

\bibitem{Lamb} P.D. Lamberti, P. Luzzini, P. Musolino, \textit{Shape perturbation of Grushin eigenvalues,}
 J. Geom. Anal., 31 (11) (2021), pp. 10679-10717.


\bibitem{Lindq} P. Linqvist, \textit{On the equation ${\rm
div}(|\nabla u|^{p-2}\nabla u)+\lambda |u|^{p-2}u=0$}, Proc. Amer.
Math. Soc. 109 (1) (1990) 157-164.

\bibitem{Lindq2} P. Lindqvist, \textit{Notes on the $p$-Laplace
equation.} Report. University of Jyv\"{a}skyl\"{a}, Department of
Mathematics and Statistics. 102 University of Jyv\"{a}skyl\"{a},
2006.


\bibitem{Liu} M. Liu, Z. Tang and C. Wang, \textit{Infinitely many solutions for a critical Grushin-type problem via local Pohozaev identities},
Ann. Mat. Pura Appl. 199 (2019), 1737-1762.

\bibitem{LiuXiao} Y. Liu, J. Xiao, \textit{Functional capacities on
the Grushin space $G_\alpha^n$}, Annali Mat. Pura Appl. (2018)
197:673702 

\bibitem{LoiuSobol} A. Loiudice, \textit{Sobolev inequalities with remainder terms for sublaplacians and other subelliptic
operators}, NoDEA Nonlinear Differential Equations Appl. 13
(2006), no. 1, 119-136.

\bibitem{Loiu00} A. Loiudice, \textit{Semilinear subelliptic problems with critical growth on Carnot groups}, Manuscripta Math. 124
(2007), 247-259.

\bibitem{LoiuNonlinear} A. Loiudice, \textit{Asymptotic behaviour of solutions for a class of degenerate
elliptic critical problems}, Nonlinear Anal. 70, no. 8 (2009),
2986-2991.

\bibitem{LoiuOpt} A. Loiudice, \textit{Optimal decay of
$p$-Sobolev extremals on Carnot groups}. J. Math. Anal. Appl. 470
no. 1, 619-631 (2019)


\bibitem{LoiuAMPA} A. Loiudice, \textit{Asymptotic estimates and nonexistence results for critical problems with
Hardy term involving Grushin-type operators}, Ann. Mat. Pura Appl.
198 (2019), 1909-1930.

\bibitem{LoiuTrends} A. Loiudice, \textit{Existence Results for Critical Problems Involving p-Sub-Laplacians on Carnot Groups}, Trends
Math. 2022, pp. 135-151.

\bibitem{LT} D.T. Luyen, N.M. Tri \textit{Existence of
 infinitely many solutions for semilinear degenerate Schr\"{o}dinger equations.}
 J. Math. Anal. Appl. 2018;461 no.2:1271-1286.

 \bibitem{MMS1} P. Malanchini P, G. Molica Bisci, S. Secchi \,\textit{A note on nonlinear
critical problems involving the Grushin subelliptic operator:
bifurcation and multiplicity results}, Potential Analysis 2025
https://doi.org/10.1007/s11118-025-10199-z

\bibitem{MMS} P. Malanchini, G. Molica Bisci, S. Secchi, \textit{Bifurcation and multiplicity results for critical problems involving the
p-Grushin operator}, Advances in Nonlinear Analysis 14(1) 2025
DOI: 10.1515/anona-2025-0089.

\bibitem{Mon} R. Monti, \textit{Sobolev inequalities for weighted gradients}, Commun. Partial Differ. Equ. 31(2006), 1479-1504.


\bibitem{Mon1} R. Monti and D. Morbidelli, \textit{Kelvin transform for Grushin operators and critical semilinear equations,}
Duke Math. J. 131 (2006), 167-202.


\bibitem{PSY} K. Perera, M. Squassina, Y. Yang, \textit{Bifurcation
and multiplicity results for critical fractional p-Laplacian
problems}, Math. Nachr. 289, No. 2-3, 332-342 (2016)  

\bibitem{R} P.H. Rabinowitz, \textit{Minimax methods in critical
point theory with application to differential equations}, CBMS
Reg. Conf. Series Math. 65, Amer. Math. Soc., Providence, R.T.
1986.

\bibitem{Serrin} J. Serrin, \textit{Local behavior of solutions
of quasi-linear equations}, Acta Math. 111 (1964), n.1, 247-302.

 \bibitem{Tri1} P. T. Thuy, N. M. Tri,
\textit{Nontrivial solutions to boundary value problems for
semilinear strongly degenerate elliptic differential equations},
NoDEA, Nonlinear Differential Equations Appl. \textbf{19} (2012),
279-298.


\bibitem{Tri} N. M. Tri, \textit{On the Grushin equation}, Mat. Zametki \textbf{63}(1)
(1998), 95-105.


\bibitem{V} D. Vassilev, \textit{$L^p$-estimates and asymptotic behavior
for finite energy solutions of extremals to Hardy-Sobolev
inequalities}, Trans. Amer. Math. Soc. 363 (2011), no. 1, 37-62.

 \bibitem{Vetois} J. V\'{e}tois, \textit{A priori estimates and application
to the symmetry of solutions for critical p-Laplace equations}, J.
Differential Equations 260, no. 1 (2016), 149-161.

\bibitem{Wang} C. Wang, Q. Wang, J. Yang, \textit{On the Grushin critical problem with a cylindrical symmetry,} Adv. Differ. Equ.
 20 (2015), 77-116.

\bibitem{Wei} Y. Wei, C. Chen, Q. Chen, H. Yang, \textit{Liouville-type theorem
for nonlinear elliptic equations involving p-Laplace-type Grushin
operators}, Math. Methods Appl. Sci. 43(1) (2020), 320-333.

\bibitem{Wil} M. Willem, \textit{Minimax Theorems,} Birkh\"{a}user, Boston 123, 1996.

\bibitem{Wil1}  M. Willem, \textit{Analyse harmonique r\'{e}elle}, Hermann, Paris,
1995.

\bibitem{Xiang} C.L. Xiang, \textit{Asymptotic behaviors of solutions to
quasilinear elliptic equations with critical Sobolev growth and
Hardy potential}, J. Differential Equations 259 (2015), 3929-3954.




\end{thebibliography}
\end{document}